\documentclass{article}

\usepackage{amsmath,amssymb}
\usepackage{mathtools}
\usepackage{mathrsfs} 
\usepackage{enumitem}
\usepackage{tikz-cd}

\usepackage{microtype}

\usepackage{amsthm}

\usepackage{hyperref}
\usepackage{cleveref}

\newtheorem{theorem}{Theorem}[section]

\newtheorem{lemma}[theorem]{Lemma}
\newtheorem{prop}[theorem]{Proposition}
\newtheorem*{prop*}{Proposition}
\newtheorem*{theorem*}{Theorem}

\theoremstyle{definition}

\newtheorem{example}[theorem]{Example}
\newtheorem{remark}[theorem]{Remark}
\newtheorem{setup}[theorem]{}

\newcommand{\Emb}{\mathrm{Emb}}

\newcommand{\Diff}{\mathrm{Diff}}
\newcommand{\Gau}{\mathrm{Gau}}
\newcommand{\Aut}{\mathrm{Aut}}

\newcommand{\R}{\mathbb R}
\newcommand{\C}{\mathbb C}
\newcommand{\N}{\mathbb N}

\newcommand{\Z}{\mathbb Z}
\newcommand{\cO}{\mathcal{O}}
\newcommand{\cP}{\mathcal{P}}
\newcommand{\cG}{\mathcal{G}}

\newcommand{\id}{\mathrm{id}}
\newcommand{\ev}{\mathrm{ev}}
\newcommand{\evol}{\mathrm{evol}}

\newcommand{\Lf}{\mathbf{L}}

\newcommand{\one}{\mathbf{1}}

\newcommand{\piorb}{\pi_1^{\mathrm{orb}}}

\newcommand{\w}[1]{\widetilde{#1}}

\newcommand{\Fr}{\mathrm{Fr}}
\DeclareMathOperator{\GL}{GL}
\DeclareMathOperator{\mO}{O}
\DeclareMathOperator{\mSO}{SO}
\DeclareMathOperator{\mU}{U}
\DeclareMathOperator{\mSp}{Sp}

\title{Orbifold diffeomorphism groups, frame bundles, and the solder form}
\author{Thomas O. Rot\footnote{CTA${}^2$, Department of mathematics, Vrije Universiteit Amsterdam, \href{mailto:t.o.rot@vu.nl}{t.o.rot@vu.nl}}, Alexander Schmeding\footnote{NTNU Trondheim, \href{mailto:alexander.schmeding@ntnu.no}{alexander.schmeding@ntnu.no}}}
\date{\today}

\begin{document}

\maketitle

\begin{abstract}

We study the diffeomorphism group of an effective orbifold. The frame bundle of an effective orbifold is a manifold and we relate the diffeomorphism group of the orbifold to a subgroup of the diffeomorphism group of its frame bundle. We determine the homotopy type of several diffeomorphism groups of orbifolds. Finally, we show that there is an infinite-dimensional principal bundle relating the diffeomorphism group of a compact orbifold to the automorphism group of its orbifold frame bundle. The resulting short exact sequence splits via frame lifting. This turns the automorphisms of the frame bundle into a semidirect product of orbifold diffeomorphisms and elements in the gauge group. 
\end{abstract}

\textbf{Keywords:} orbifolds, diffeomorphism groups, frame bundles, homotopy type, solder forms, equivariant diffeomorphisms\\
\textbf{MSC2020:} 57R18 (primary); 58D05, 55P15, 22E65

\section{Introduction}

Diffeomorphism groups of manifolds encode important geometric and topological information about the underlying manifold and there is a vast literature devoted to the study of algebraic and topological features of these groups. Orbifold diffeomorphism groups are relatively less studied, but see~\cite{borzbruns,S15,pasquottorot1}. The purpose of this article is to develop two complementary descriptions of orbifold diffeomorphism groups and to use them for homotopy-theoretic and Lie-theoretic calculations.

Our first description uses the (orbifold) frame bundle. If $\cO$ is an $n$-dimensional effective orbifold, its frame bundle $\Fr(\cO)$ is a manifold which carries a locally free and proper $\GL(n)$-action. As is well known, the quotient of the frame bundle by this action is diffeomorphic to the orbifold, cf.~\cite{moerdijkmrcun,HaM04}. A $\GL(n)$-equivariant diffeomorphism of the frame bundle is called an automorphism of the frame bundle. It induces a diffeomorphism of the underlying orbifold which gives the exact sequence
\begin{equation}
\label{eq:mainseq}
 \mathbf{1}\rightarrow \Gau(\Fr(\cO))\rightarrow \Aut(\Fr(\cO))\rightarrow \Diff(\cO)\rightarrow \mathbf{1}.
\end{equation}
Here the gauge group consists  of those automorphisms that descend to the identity. This exact sequence is split: Any orbifold diffeomorphism $f\colon \cO\rightarrow \cO$ lifts to an automorphism $\Fr(f)\colon \Fr(\cO)\rightarrow \Fr(\cO)$ by its action on the frames, cf.~Proposition~\ref{prop:orb1}. The frame bundle of the orbifold carries a tautological form $\theta$, the solder form, which was studied by Daza~\cite{daza} in the context of $G$-structures on orbifolds. Theorem~\ref{thm:Fr_gp_mor} states that the group of automorphisms that preserve this solder form is the image of $\Fr$. This is the orbifold analogue of the characterisation of diffeomorphisms as frame-bundle automorphisms preserving the solder form, cf.~Theorem~\ref{thm:splitting_bundle_mfd}.

Our second description uses the universal cover of a good orbifold. For a connected good orbifold $\cO$, i.e.~an orbifold whose universal cover is a manifold $M$ with group of deck transformations $G$, we can relate the diffeomorphism groups of $\cO$ and $M$. We show that there is an exact sequence
\[
\mathbf{1}\rightarrow G\rightarrow N_{\Diff(M)}(G)\rightarrow \Diff(\cO)\rightarrow \mathbf{1}
\]
where $N_{\Diff(M)}(G)$ denotes the normalizer of $G$ in $\Diff(M)$. This generalizes an observation made in~\cite{S14}. Unlike the frame-bundle sequence, the sequence  need not split. Both perspectives on the orbifold diffeomorphism groups are naturally related, see Theorem~\ref{thm:twoperspectives}.

We then investigate homotopy types of diffeomorphism groups of several orbifolds. The (manifold) diffeomorphism group $\Diff(\R^n)$ is homotopy equivalent to $ \mO(n)$ via the Alexander trick. We generalize this calculation to orbifolds and determine the homotopy type of the orbifold diffeomorphism group of orbifolds of the form $\R^n/G$, for a finite group $G$, in terms of representation-theoretic data of $G$ in Theorem~\ref{thm:diffeochart}. 

Martin Pinsonnault asked on MathOverflow~\cite{pinsonnault} what is known about the homotopy type of diffeomorphism groups of two-dimensional orbifolds. In~\cite{earlemcmullen} and~\cite{earle} it is shown that the diffeomorphism groups of hyperbolic two-dimensional orbifolds have contractible components. As far as we are aware, not much is known for the bad two-dimensional orbifolds. We show in Propositions~\ref{prop:teardrop} and~\ref{prop:spindle} that the homotopy type of the diffeomorphism groups of the teardrop and spindle orbifolds is $\mO(2)$. 

The final part of the article concerns infinite-dimensional Lie group structures. The orbifold diffeomorphism group $\Diff(\cO)$ is an infinite-dimensional Lie group,~\cite{S15}. We equip $\Aut(\Fr(\cO))$ and $\Gau(\Fr(\cO))$, for compact $\cO$, with a regular Fr\'echet Lie group structure such that Sequence~\eqref{eq:mainseq} is a split exact sequence of regular Fr\'echet Lie groups. In particular $\Aut(\Fr(\cO))$ is a nontrivial semidirect product and~\eqref{eq:mainseq} describes a trivial infinite-dimensional principal bundle.  

Note that $\Diff(\Fr(\cO))$ is also an infinite-dimensional Lie group, in a fine Whitney-type topology (since $\Fr(\cO)$ is non-compact, cf.~\cite{Mic80}). This topology is too fine for our purpose as it makes the subgroup identified with $\Diff(\cO)$ via frame lifting into a discrete topological group. Hence, it is important to note that the Lie group structure of the semidirect product group $\Aut(\Fr(\cO))$ is not inherited from $\Diff(\Fr(\cO))$. This phenomenon is well known in the manifold world. If $\cO$ is an ordinary compact manifold, our result specializes to the automorphism group of its principal frame bundle and recovers the classical principal-bundle extension studied in~\cite{AaCaMaM89}. See Remark~\ref{rem:old_bundle} for more information.

\subsection*{AI statement}

We began this project a year before finalizing this paper. The original goal, conceived before we interacted with LLMs, was to find an orbifold analogue of the classical result of Theorem~\ref{thm:splitting_bundle_mfd}, and to use this to study the homotopy types of orbifold diffeomorphism groups and answer the MathOverflow question of Martin Pinsonnault~\cite{pinsonnault}. 

Due to other obligations, we did not work on this project until the summer of 2026. With the state of mathematics in flux we felt the urge to experiment with LLMs. Many, if not most, of the proofs of this paper were therefore conceived through interactions with various LLMs. None of the text in the paper is generated by LLMs. 

Thomas Rot declares that he has received free access to state-of-the-art model(s) during the editing phase of this paper. 

It does not need to be stated that the authors take full responsibility for the results. 

\subsection*{Acknowledgment}
Thomas Rot would like to thank Bob Rink for a useful discussion on the isotypical decomposition of a representation. 

\section{The manifold frame bundle and diffeomorphisms}

We recall some standard material on the frame bundle and the solder form of a manifold. See for example~\cite[I.5 Example 5.1 and II.2]{KaN63}.

Let $M$ be an $n$-dimensional manifold and $\pi\colon\Fr(M)\rightarrow M$ its frame bundle. This is the principal $\GL(n)$-bundle whose fiber over $x\in M$ is
\[
\Fr_x(M)=\{p\colon \R^n\rightarrow T_xM\,|\,p\text{ is a linear isomorphism}\}.
\]
 An element in $\Fr(M)$ is denoted by $p_x=(x,p)$ where $p\in \Fr_x(M)$. The right action of $\GL(n)$ is given by the formula $p_xA=(x,pA)=(pA)_x$. We denote by $\nu_A\colon\Fr(M)\rightarrow \Fr(M)$ the map $p_x\mapsto p_xA$.
 The \emph{solder form} is the $\R^n$-valued one-form $\theta\in \Omega^1(\Fr(M),\R^n)$ defined by the equation $\theta_{p_x}=p^{-1}T_{p_x}\pi$. If we need to emphasize the underlying manifold $M$ we denote the solder form by $\theta^M$. We recall some basic facts.

\begin{prop}\label{prop:solder0}
 Let $p_x\in \Fr(M)$ and $A\in \GL(n)$. The solder form satisfies
 \begin{enumerate}
 \item[i)] $\nu_A^*\theta=A^{-1}\theta$.
 \item[ii)] $\ker \theta_{p_x}=T_{p_x}(p_x\GL(n))$.
 \item[iii)] $\theta_{p_x}\colon T_{p_x}\Fr(M)\rightarrow \R^n$ is surjective. 
 \end{enumerate}
\end{prop}
\begin{proof}
 As $\pi\nu_A=\pi$ we have that $T\pi T\nu_A=T\pi$. We then compute
 \begin{align*} 
 \nu_A^*(\theta)_{p_x}=\theta_{p_xA}T_{p_x}\nu_A=A^{-1}p^{-1}T_{p_x A}\pi \,T_{p_x}\nu_A =A^{-1}p^{-1}\,T_{p_x}\pi=A^{-1}\theta_{p_x}
 \end{align*}
 which proves i). Since $p$ is a linear isomorphism we find that $\ker \theta_{p_x}=\ker T_{p_x}\pi$, which is exactly the tangent space to the orbit $T_{p_x}(p_x\, \GL(n))$. Finally $\theta_{p_x}$ is surjective as it is the composition of $T_{p_x}\pi$ and $p^{-1}$, which are both surjective. 
\end{proof}

\begin{prop}\label{prop:solder1}
Let $N$ be another manifold and let $f\colon M\rightarrow N$ be a local diffeomorphism. Then $f$ induces a local diffeomorphism $\Fr(f)\colon\Fr(M)\rightarrow \Fr(N)$ defined by $\Fr (f)p_x=(f(x),T_xf\,p)$. The assignment $f\mapsto \Fr(f)$ is functorial, and if $f$ is a diffeomorphism, so is $\Fr(f)$. Moreover 
\begin{enumerate}
 \item[i)] $\Fr(f)$ is equivariant with respect to the $\GL(n)$-action.
 \item[ii)] $\Fr(f)$ preserves the solder form: $\Fr(f)^*\theta^N=\theta^M.$
\end{enumerate}

\end{prop}
\begin{proof}
Note that $\Fr$ is functorial by the chain rule. If $f$ is a diffeomorphism then $\Fr(f)\circ \Fr(f^{-1})=\Fr(\id_N)=\id_{\Fr(N)}$ shows that $\Fr(f)$ is a diffeomorphism. If $f$ is only a local diffeomorphism, it restricts to a diffeomorphism between suitable open sets. Hence $\Fr(f)$ is always a local diffeomorphism.

Equivariance follows from the calculation 
\begin{align*}
\Fr(f)(p_xA)=(f(x),T_xf\,p\, A)=(f(x),T_xf\, p) A
=(\Fr(f)(p_x)) A.
\end{align*}
Note that $f\circ \pi^M=\pi^N\circ\Fr(f)$ so $Tf\, T\pi^M=T\pi^N\, T \Fr(f)$, and therefore
\begin{align*}
\left(\Fr(f)^*\theta^N\right)_{p_x}&=(T_xf\,p)^{-1}T_{\Fr(f)(p_x)}\pi^N \; T_{p_x}\Fr(f)\\
&=p^{-1}(T_xf)^{-1}T_xf\, T_{p_x}\pi^M=\theta^M_{p_x}.\qedhere
\end{align*}
\end{proof}
Conversely we have
\begin{prop}\label{prop:solder2}
 Let $F\colon\Fr(M)\rightarrow \Fr(N)$ be an equivariant local diffeomorphism. Then $F$ induces a unique local diffeomorphism $f\colon M\rightarrow N$. The association $F\mapsto f$ is functorial. Moreover
 \begin{enumerate}
 \item[i)] If $F$ is a diffeomorphism, then $f$ is as well.
 \item[ii)] If $F$ preserves the solder form $F^*\theta^N=\theta^M$, then $F=\Fr(f)$. 
 \end{enumerate}
\end{prop}

\begin{proof}
Because $F$ is equivariant it descends to a smooth map $f\colon M\rightarrow N$ such that $\pi^N F=f\pi^M$. As $\pi^M$ is surjective, $f$ is uniquely determined. Now $T\pi^N TF=Tf T\pi^M$. As $T\pi^N$ and $TF$ are surjective, it follows that $Tf$ is surjective, hence invertible. Thus $f$ is a local diffeomorphism if $F$ is. Clearly, if $F'$ is another equivariant local diffeomorphism that descends to $f'$, then $F'\circ F$ descends to $f'\circ f$. If $F$ is a diffeomorphism, then $F^{-1}F$ descends to the identity, thus $f$ is a diffeomorphism if $F$ is. 

Now assume that $F$ preserves the solder form. For $p_x\in \Fr(M)$ we write $F(p_x)=(f(x),\widehat F(p_x))$. We then compute
\begin{align*}
(F^*\theta^N)_{p_x}&=(\widehat F(p_x))^{-1}T_{F(p_x)}\pi^N\; T_{p_x}F=(\widehat F(p_x))^{-1}T_xf\; T_{p_x}\pi^M\\
&=(\widehat F(p_x))^{-1}T_xf\, p \,\theta^M_{p_x}.
\end{align*}
As $\theta^M_{p_x}$ is surjective, cf.~Proposition~\ref{prop:solder0}, the right-hand side is equal to $\theta^M_{p_x}$ if and only if $(\widehat F(p_x))^{-1} T_xf \,p$ is the identity. Thus $\widehat F(p_x)=T_xf\,p$ and $F=\Fr(f).$ 
\end{proof}

For a smooth manifold $M$ we consider 
\begin{align*}
 \Aut(\Fr(M))&=\{F\in \Diff(\Fr(M))\,|\,F(p_x A)=F(p_x) A\\ & \qquad \text{ for all } A\in \GL(n) \text{ and } p_x\in \Fr(M)\}
\end{align*} 
the \emph{automorphism group}. These are all diffeomorphisms that are $\GL(n)$-equivariant. Further, denote by 
\begin{align*}
\Gau(\Fr(M))=\{F\in \Aut(\Fr(M))\,|\,\pi F=\pi\}\end{align*} the \emph{gauge group}, i.e.~those automorphisms that descend to the identity. 
Propositions~\ref{prop:solder1} and~\ref{prop:solder2} show that we have an exact sequence of groups
\begin{align}\label{eq:gauge_ex_seq}
\mathbf{1}\rightarrow \Gau(\Fr(M))\rightarrow \Aut(\Fr(M))\rightarrow \Diff(M)\rightarrow \mathbf{1}
\end{align}
that is split on the right by the map $f\mapsto \Fr(f)$. The image of the splitting is
\begin{equation*}
 \Aut(\Fr(M),\theta)=\{F\in \Aut(\Fr(M))\,|\,F^*\theta=\theta\},
\end{equation*}
those automorphisms that preserve the solder form. Thus we obtain
\begin{theorem}\label{thm:splitting_bundle_mfd}
The map $\Fr \colon \Diff(M) \rightarrow \Diff(\Fr(M)), f\mapsto \Fr(f)$ induces an isomorphism between $\Diff(M)$ and $\Aut(\Fr(M),\theta)$. 
\end{theorem}
The sequence \eqref{eq:gauge_ex_seq} is a special case of the results obtained in~\cite{AaCaMaM89} with Theorem~\ref{thm:splitting_bundle_mfd} providing a splitting section for this sequence, cf.~Remark~\ref{rem:old_bundle}. We shall later lift the splitting from the abstract group level to a splitting of (infinite-dimensional) Lie groups.

\section{The orbifold frame bundle and diffeomorphisms}\label{sect:ofd_frame bundle}

We now generalize the results of the previous section to orbifolds. All orbifolds we will consider are effective and we drop the adjective effective in what follows. Let $\cO$ be an $n$-dimensional orbifold and $O$ its underlying space. To fix notation we recall that an \emph{orbifold chart} on $U\subseteq O$ is a triple $(\w U, \Gamma,\varphi)$, where $\w U$ is an $n$-dimensional manifold, $\Gamma$ is a finite group acting smoothly and effectively on each component of $\w U$ and $\varphi\colon \w U\rightarrow U\subseteq O$ a map inducing a homeomorphism $\w U/\Gamma\rightarrow U$. A \emph{change of charts} between $(\w U_i,\Gamma_i,\varphi_i)$ and $(\w U_j,\Gamma_j,\varphi_j)$ is a smooth embedding $\lambda\colon \w U_i\rightarrow \w U_j$ such that $\varphi_j \lambda=\varphi_i$. Two charts are \emph{compatible at} $x\in U_i\cap U_j$, if there exists a third chart $(\w U_k,\Gamma_k,\varphi_k)$ with $x\in U_k\subseteq U_i\cap U_j$ that embeds into both $\w U_i$ and $\w U_j$. An \emph{orbifold atlas} is a compatible collection of orbifold charts that covers $O$. One atlas $\mathcal{A}$ \emph{refines} the atlas $\mathcal B$ if every chart in $\mathcal A$ embeds into some chart in $\mathcal B$. Two atlases are \emph{equivalent} if there exists a common refinement, and just as in the manifold case, an orbifold admits a maximal atlas.

Let $x\in \cO$ and $U$ a neighborhood with chart $(\w U,\Gamma,\varphi)$ and $\w x\in \w U$ such that $\varphi(\w x)=x$. The \emph{isotropy group} at $\w x$ is defined to be $\Gamma_{\w x}=\{g\in \Gamma\,|\, g\w x=\w x\}$. Up to isomorphism, this group is independent of the chart and lift $\w x$. We will therefore abuse notation and denote this by $\Gamma_x$ when no confusion can occur. If $\Gamma_x=1$ is the trivial group, then we call $x$ a \emph{smooth point}. 
\subsection{The orbifold frame bundle and solder form}
 
Let $(\w U,\Gamma, \varphi)$ be a chart. For every $g\in \Gamma$, the map $\mu_g\colon \w U\rightarrow \w U$, $\mu_g\w x=g\w x$ is a diffeomorphism, hence induces a map $\Fr(\mu_g)\colon \Fr(\w U)\rightarrow \Fr(\w U)$. Thus $\Gamma$ acts on $\Fr(\w U)$ by the formula $g\,p_x=\Fr(\mu_g)(p_x).$ 
 
\begin{lemma}
\label{lem:fract}
 Let $(\w U, \Gamma, \varphi)$ be an orbifold chart. The action of $\Gamma$ is free on $\Fr(\w U)$ and commutes with the $\GL(n)$-action. 
\end{lemma}
\begin{proof}
Without loss of generality, we may assume $\w U$ is connected. Choose an auxiliary invariant Riemannian metric on $\w U$. Then $\Gamma$ acts by isometries on $\w U$. Now suppose that $g\, {\w p}_{\w x}={\w p}_{\w x}$ in $\Fr(\w U)$. Then $g\, \w x=\w x$ and $T_{\w x}\mu_g\w p=\w p$. As $\w p$ is a frame, this implies that the differential $T_{\w x}\mu_g$ is the identity on $T_{\w x}\w U$. An isometry of a connected manifold that fixes a point and whose differential at this point is the identity, is the identity map, cf.~\cite[Proposition 5.6.2]{Pet16}. Since the action is assumed to be effective, the action is free. 
 
 Now let $g\in \Gamma$ and $A\in \GL(n)$. Then $g \cdot ({\w p}_{\w x}\cdot A)=\Fr(\mu_g)(\w x,\w p A)=(\Fr(\mu_g)({\w p}_{\w x}))\cdot A=(g\cdot {\w p}_{\w x})\cdot A$, thus the actions of $\Gamma$ and $\GL(n)$ commute.
\end{proof}
 
The action of $\Gamma$ is free, hence $\Fr(U):=\Fr(\w U)/\Gamma$ is a manifold. We denote the (covering) projection by $\Fr(\varphi)\colon \Fr(\w U)\rightarrow \Fr(U)$ and note that it is a submersion. As the actions of $\Gamma$ and $\GL(n)$ commute, the projection $\Fr(\varphi)$ is $\GL(n)$-equivariant. A change of charts induces an embedding $\Fr(U_i)\rightarrow \Fr(U_j)$, which can be used to form the frame bundle of the orbifold $\pi\colon \Fr(\cO)\rightarrow \cO$. We emphasize that $\Fr(\cO)$ is a manifold. As the actions of the local groups commute with the $\GL(n)$-action, the group $\GL(n)$ acts on $\Fr(\cO)$ from the right. We denote by $\nu_A\colon \Fr(\cO)\rightarrow \Fr(\cO)$ the map $\nu_A(p_x)=p_x A.$ This action is not free (unless $\cO$ happens to be a manifold). The stabilizer of $p_x\in \Fr(\cO)$ is $\GL(n)_{p_x}=\{A\in \GL(n)\,|\,p_x A=p_x\}$ and is isomorphic to $\Gamma_{x}.$ Thus the action is locally free.

The action is also proper and the (orbifold) quotient is well known to be diffeomorphic to $\cO.$ See~\cite[Section 4.2]{moerdijkmrcun} or~\cite{HaM04} for more details on this construction. The orbifold frame bundle carries a tautological form which we will call the \emph{solder form}. 
We remark that the solder form on the frame bundle of an orbifold was defined and studied by Daza~\cite{daza} in the context of $G$-structures on orbifolds. 
 
\begin{prop}\label{prop:solder_orbi}
 There exists a canonical \emph{solder form} $\theta\in \Omega^1(\Fr(\cO),\R^n)$ on $\Fr(\cO)$ that satisfies
 \begin{enumerate}
 \item[i)] $\nu_A^*\theta=A^{-1}\theta$
 \item[ii)] $\ker \theta_{
 p_x}=T_{p_x}\left(p_x\GL(n)\right)$
 \item[iii)] $\theta_{p_x}$ is surjective.
 \end{enumerate} 
\end{prop}
\begin{proof}
 Denote the solder form on $\Fr(\w U)$ by $\theta^{\w U}\in \Omega^1(\Fr(\w U),\R^n)$. By Proposition~\ref{prop:solder1} the solder form is invariant under the action of $\Gamma$, i.e.~$\Fr(\mu_g)^*\theta^{\w U}=\theta^{\w U}$. It follows that the solder form descends to an $\R^n$-valued form $\theta^U$ on $\Fr(U)$. This form is independent of the chosen orbifold chart of $U$: Indeed, if $\lambda\colon \w U_i\rightarrow \w U_j$ is a change of charts, $\Fr(\lambda)^*\theta^{\w U_j}=\theta^{\w U_i}$ by Proposition~\ref{prop:solder1}. This implies that the local forms $\theta^U$ patch together to a global form $\theta$ on $\Fr(\cO)$. Properties i), ii) and iii) follow directly from the corresponding properties of $\theta^{\w U}$ in Proposition~\ref{prop:solder0}.
\end{proof}

\subsection{Orbifold diffeomorphisms and the frame bundle}
Orbifold morphisms are continuous maps of the underlying topological spaces which lift smoothly in charts and interact suitably with the isotropy information. Even for effective orbifolds general mappings are technical and there is a variety of different approaches to them, cf.~\cite{BaB13,Car19,Poh17}. Fortunately, the different notions tend to yield the same orbifold diffeomorphisms for effective orbifolds. Therefore, we do not dwell on these different concepts here but see~\cite{S15} for a discussion.

Let $\cO$ and $\cP$ be orbifolds with underlying spaces $O$ and $P$. A continuous map $f\colon O\rightarrow P$ defines a \emph{smooth orbifold map} if for every $x\in O$ there exist neighborhoods $U\subseteq O$ and $V\subseteq P$ containing $x$ and $f(x)$ respectively, and corresponding orbifold charts $(\w U,\Gamma,\varphi)$ and $(\w V,\Lambda, \psi)$ to which $f$ admits a smooth local lift $\widetilde f\colon \w U\rightarrow \w V$ such that $f\varphi=\psi\w f $. We will write $f\colon \cO\rightarrow \cP$ for a smooth orbifold map. The map $f$ is an \emph{orbifold diffeomorphism} if there exists a smooth orbifold map $f^{-1}\colon \cP\rightarrow \cO$ such that $f^{-1}\circ f=\id_\cO$ and $f\circ f^{-1}=\id_\cP$. 

A local lift of an orbifold diffeomorphism is a local diffeomorphism~\cite[Lemma 2]{pasquottorot1} and, by shrinking the charts as necessary, the local lift is a diffeomorphism. We recall the following behaviour of the local lift~\cite[Lemma 3]{pasquottorot1}. If $f$ is an orbifold diffeomorphism, and $\w{f}$ is a smooth diffeomorphism lifting $f$ between the connected orbifold charts $(\w{U},\Gamma,\varphi)$ and $(\w{V},\Lambda,\psi)$, then there is an induced group isomorphism $\alpha\colon \Gamma\rightarrow \Lambda$ such that $\w{f}(g\, \w{x})=\alpha(g)\, \w{f}(\w{x})$ for all $\w{x}\in \w{U}$ (in other words, $\w{f}$ is $\alpha$-equivariant). If $\w f'$ is another lift of $f$ between the same orbifold charts, then $\w f=g\, \w f'$ for a unique $g\in \Lambda$. 

Analogously to the manifold case, we define the following subgroups of the diffeomorphism group of the orbifold frame bundle: 
\begin{align*}
\Aut(\Fr(\cO)) &= \{F \in \Diff(\Fr(\cO))\,|\, F(p_x A)=F(p_x) A, \\&\qquad \text{ for all } p_x\in \Fr(\cO),\text{ and } A \in \GL(n)\}\\
\Aut(\Fr(\cO),\theta)&=\{F\in \Aut(\Fr(\cO))\,|\,F^*\theta=\theta\}
\end{align*} 
the group of \emph{automorphisms} and those that preserve the solder form.

We now show that all orbifold diffeomorphisms lift to automorphisms of the frame bundle. This was also observed by Kankaanrinta~\cite[Lemma 3.2]{Kankaanrinta}.
\begin{prop}
\label{prop:orb1}
 Let $f\in \Diff(\cO)$. Then $f$ induces an automorphism $\Fr(f)\in \Aut(\Fr(\cO))$ such that $\Fr(f)^*\theta=\theta$. 
\end{prop}
\begin{proof}
Let $\w f\colon \w U\rightarrow \w V$ be a local lift to connected charts $(\w U,\Gamma,\varphi)$, $(\w V,\Lambda,\psi)$. Then $\w f$ is a local diffeomorphism and thus induces a map $\Fr(\w f)\colon \Fr(\w U)\rightarrow \Fr(\w V)$. The map is $\alpha$-equivariant hence induces a map
$\Fr(f)\bigr|_{\Fr(U)}\colon \Fr(U)\rightarrow \Fr(V)$. Note that this map does not depend on the choice of lift between $\w U$ and $\w V$. If $\w f'$ is another local lift of $f$ to the same charts, then $\w f'=g\, \w f$ for a unique $g\in \Lambda$. But then $\Fr(\w f')=\Fr(g\, \w f)=\Fr(\mu_g)\Fr(\w f)=g\, \Fr(\w f)$ and it follows that $\Fr(f)\bigr|_{\Fr(U)}$ does not depend on the lift $\w f$. A similar argument, applied to change of charts, shows that we can patch $\Fr(f)\bigr|_{\Fr(U)}$ to a smooth map $\Fr(f)\colon \Fr(\cO)\rightarrow \Fr(\cO).$ Exploiting functoriality we get that
\[\id_{\Fr(\cO)}=\Fr(\id_\cO)=\Fr(f\circ f^{-1})=\Fr(f)\circ \Fr(f^{-1})\] and it follows that $\Fr(f)$ is a diffeomorphism. As $\Fr(\widetilde f)$ is $\GL(n)$-equivariant, so is $\Fr(f)$. Since $\Fr(\w f)^*\theta^{\w V}=\theta^{\w U}$, by Proposition~\ref{prop:solder1}, it follows that $\Fr(f)^*\theta=\theta$.
\end{proof}

\begin{prop}
\label{prop:orb2}
 Let $F\in \Aut(\Fr(\cO))$. Then $F$ descends to a unique orbifold diffeomorphism $f\in \Diff(\cO).$ If $F$ preserves the solder form~$F^*\theta=\theta$, then $F=\Fr(f)$.
\end{prop}
\begin{proof}
As $F$ is $\GL(n)$-equivariant, it descends to a continuous map $f\colon O\rightarrow O$ such that $\pi F=f\pi$. As $\pi$ is surjective, $f$ is uniquely determined. We claim that $f$ is an orbifold diffeomorphism.
 
A first observation is that $F^{-1}F=\id=F\circ F^{-1}$, and since the association $F\mapsto f$ is functorial, it follows that $f\colon \cO\rightarrow \cO$ admits an inverse $f^{-1}$. We now show that $f$ (and $f^{-1}$) is a smooth orbifold map. 

Let $x\in O$ and choose a chart $(\w V,\Lambda,\psi)$ around $f(x)$ and a linear chart $(\w U=\R^n,\Gamma,\varphi)$ around $x$ such that $f(U)\subseteq V.$ As $\w U= \R^n$, the frame bundle $\Fr(\w U)\rightarrow \w U$ is trivial. Let $s\colon \w U\rightarrow \Fr(\w U)$ be a smooth section. As $\Fr(\psi)\colon \Fr(\w V)\rightarrow \Fr(V)$ is a smooth covering space, and $\w U$ is simply connected, there exists a lift $\sigma$, as in the diagram below
\[
\begin{tikzcd}
&& &\Fr(\w V)\arrow[d, "\Fr(\psi)"]\\
\w U \arrow[rrru,dashed,"\sigma"]\arrow[rrr, "F\circ \Fr(\varphi)\circ s"']&&&\Fr(V)
\end{tikzcd}
\]
which is readily seen to be smooth, cf.~\cite[Proposition 4.33]{Lee13}. Let $\w f\colon \w U\rightarrow \w V$ be $\pi_{\w V}\circ \sigma$ where $\pi_{\w V}\colon \Fr(\w V)\rightarrow \w V.$ Then $\w f$ is smooth and
\[
\psi \circ \w f=\psi \circ \pi_{\w V} \circ \sigma= \pi_{V}\circ \Fr(\psi) \circ \sigma=\pi_V \circ F\circ \Fr(\varphi)\circ s=f\circ \pi_U\circ \Fr(\varphi)\circ s=f\circ \varphi.
\]
The map $\w f$ is a smooth lift of $f$. We conclude that $f$, and similarly $f^{-1}$, is a smooth orbifold map, hence an orbifold diffeomorphism. 

We now define $\w F\colon \Fr(\w U)\rightarrow \Fr(\w V)$ by $\w F(s(\w x)\, A)=\sigma (\w x)\, A.$ Then $\w F$ is $\GL(n)$-equivariant by construction and $\w F$ lifts $F$. Indeed
\[
\Fr(\psi)\circ \w F(s(\w x)\, A)=\Fr(\psi)(\sigma(\w x)\, A) = F\circ \Fr(\varphi)(s(\w x)\, A).
\]
The induced map $\w U\rightarrow \w V$ is $\w f$. If $F$ preserves $\theta$ then
\[
\w F^*\theta^{\w V}=\w F^*\Fr(\psi)^*\theta=\Fr(\varphi)^*F^*\theta=\Fr(\varphi)^*\theta=\theta^{\w U},
\]
hence $\w F=\Fr(\w f)$ by Proposition~\ref{prop:solder2}.  Then
\[
F\circ \Fr(\varphi)=\Fr(\psi)\w F=\Fr(\psi)\Fr(\w f)=\Fr(f)\circ \Fr(\varphi).
\]
Thus $F=\Fr(f)$ as $\Fr(\varphi)$ is surjective.
\end{proof}
 The associations in Propositions~\ref{prop:orb1} and~\ref{prop:orb2} are functorial, thus we obtain the exact sequence of groups
\begin{equation}
 \label{eq:frameseq}
\mathbf{1}\rightarrow \Gau(\Fr(\cO))\rightarrow \Aut(\Fr(\cO))\rightarrow \Diff(\cO)\rightarrow \mathbf{1},
\end{equation}
where the \emph{gauge group} is
$\Gau(\Fr(\cO))=\{F\in \Aut(\Fr(\cO))\,|\,\pi F=\pi\}$.
The sequence is right split by the map $f\mapsto \Fr(f)$. The image of this map is exactly $\Aut(\Fr(\cO),\theta).$ We obtain our first main result. 
 
\begin{theorem}\label{thm:Fr_gp_mor}
 $\Fr\colon \Diff(\cO)\rightarrow \Aut(\Fr(\cO),\theta)$ is a group isomorphism. 
\end{theorem}

\section{The diffeomorphism group of a good orbifold}

Theorem~\ref{thm:Fr_gp_mor} describes the orbifold diffeomorphism group in terms of a subgroup of diffeomorphisms of the frame bundle. If the orbifold is \emph{good}, i.e.~an orbifold covered by a manifold, then the situation is simpler. The goal of this section is to explain this. We first recall some algebraic constructions. 

The \emph{normalizer} of a subgroup $G$ of a group $H$ is the subgroup \[N_H(G)=\{h\in H\,|\,hGh^{-1}=G\}.\] The normalizer is the largest subgroup of $H$ that contains $G$ as a normal subgroup. We define the map
$c\colon N_H(G)\rightarrow \mathrm{Aut}(G)$ via $c(h)=c_h$ where $c_h(g)=h g h^{-1}$. The kernel of $c$ is the \emph{centralizer} $C_H(G)=\{h\in H\,|\, hg=gh \text{ for all } g\in G\}$. Thus we have an exact sequence
\begin{equation}
\label{eq:centnorm}
\mathbf{1}\rightarrow C_H(G)\rightarrow N_H(G)\rightarrow \mathrm{Aut}(G).
\end{equation}
We can directly verify that another description of the normalizer is
\[N_H(G)=\{h\in H\,|\, \text{ there exists } c_h\in \mathrm{Aut}(G) \text{ such that }hg=c_h(g)h\,\text{ for all } g\in G\}.\] 
Now let $\cO$ be a connected good orbifold. Let $r\colon M\rightarrow \cO$ be its universal (orbifold) cover,~\cite[Chapter 13, Proposition 13.2.4]{thurston}. As the orbifold is good, $M$ is a manifold. The group of deck transformations
\[
G=\{g\in \Diff(M)\,|\, r\circ g=r\}
\]
is (depending on one's viewpoint even by definition) isomorphic to the orbifold fundamental group $\piorb(\cO)$. 

\begin{lemma}\label{lem:OFD1}
 Let $r \colon M \rightarrow \cO$ be the connected universal orbifold cover of a good orbifold and $F\in N_{\Diff(M)}(G)$. Then $F$ descends to an orbifold diffeomorphism $f\colon \cO\rightarrow \cO.$ 
\end{lemma}
\begin{proof}
An element $F\in N_{\Diff(M)}(G)$ descends to a map $f\colon O\rightarrow O$ as it maps orbits to orbits. We claim that $f$ is an orbifold diffeomorphism. Choose $x\in \cO$ and $\w x\in M$ with $r(\w x)=x.$ Let $G_{\w x}=\{g\in G\,|\, g\w x=\w x\}$ be the stabilizer. This is a finite group. Let $\w U$ be a neighborhood of $\w x$ such that $G_{\w x}\, \w U=\w U$ and $(g\, \w U)\cap\w U=\emptyset$ for all $g\not \in G_{\w x}$, cf.~\cite[Proposition 13.2.1]{thurston}. Then $(\w U, G_{\w x},r)$ and $(F(\w U), G_{F(\w x)},r)$ are orbifold charts for $\cO$, and $F$ is a smooth lift of $f$ to these charts. Applying the same argument to $f^{-1}$ shows that $f$ is an orbifold diffeomorphism. 
\end{proof}

We also have the converse, cf.~\cite[Proposition 9]{choi}.

\begin{lemma}\label{lem:OFD2}
 Let $f\in \Diff(\cO)$, $x\in \cO$ be a smooth point, and $\w x \in r^{-1}(x)$ a lift. Then given a lift $\w y\in r^{-1}(y)$ of $y=f(x)$, there is a unique lift $F\in N_{\Diff(M)}(G)$ of $f$ with $F(\w x)=\w y$.
\end{lemma}
\begin{proof}
Observe that $f\circ r\colon M\rightarrow \cO$ is a universal orbifold covering. By orbifold covering theory the universal orbifold cover is unique up to isomorphism, and therefore there is a diffeomorphism $F'\colon M\rightarrow M$ such that $r\circ F'=f\circ r$. Now $F'(\w x)\in r^{-1}(y)$. 

Since $f$ is an orbifold diffeomorphism, $y=f(x)$ is also a smooth point and the group of deck transformations $G$ acts freely and transitively on the fiber $r^{-1}(y)$. 

By transitivity we can choose $g\in G$ such that $g(F'(\w x))=\w y$. Define $F=g \circ F'$. Then $F$ is a lift of $f$ with $F(\w x)=\w y$. 

If $F''$ is another lift of $f$ with $F''(\w x)=\w y$, then $r\circ F\circ (F'')^{-1}=f\circ f^{-1}\circ r=r$. Thus $F\circ (F'')^{-1}\in G$ and $F \circ (F'')^{-1}(\w y)=\w y$. By freeness of the action of the group of deck transformations on the fiber, we find that $F$ must be the unique lift of $f$ with $F(\w x)=\w y$. 

As $r F=fr$ we have that $r\circ F\circ g\circ F^{-1}=f\circ r\circ g\circ F^{-1}=f\circ r\circ F^{-1}=f\circ f^{-1}\circ r=r$ for all $g\in G$. Thus $F g F^{-1}\in G$, i.e. $F\in N_{\Diff(M)}(G)$.
\end{proof}

Lemmas~\ref{lem:OFD1} and~\ref{lem:OFD2} combine to give the exact sequence
\begin{equation}
\label{eq:coveringsequence}
\mathbf{1}\rightarrow G\rightarrow N_{\Diff(M)}(G)\rightarrow \Diff(\cO)\rightarrow \mathbf{1}.
\end{equation}

\begin{remark}
In contrast with the frame bundle sequence \eqref{eq:frameseq}, which is split by the frame map $\Fr$, this sequence does \emph{not} split in general. This already fails in the manifold case. To see this explicitly, consider the universal cover $r\colon \R\rightarrow S^1$, $r(x)=(\cos(2\pi x),\sin(2 \pi x))$. Let $f\in \Diff(S^1)$ denote the rotation by angle $\pi$. Then $f^2=\id_{S^1}$. But for any lift $F$ of $f$, which must be of the form $F(x)=x+n+\frac12$ for some $n\in \Z$, the map $F^2$ is not the identity. 
\end{remark}
We obtain the following proposition, which generalizes~\cite[Corollary 7]{S14}.
\begin{prop}
\label{prop:universal}
 Let $\cO$ be a connected good orbifold with universal cover $r\colon M\rightarrow \cO$ and group of deck transformations $G=\piorb(\cO)$. Then $\Diff(\cO)\cong N_{\Diff(M)}(G)/G$.
\end{prop}

We now have two descriptions of the orbifold diffeomorphism group for a good orbifold, one via the frame bundle, cf.~Theorem~\ref{thm:Fr_gp_mor} and one via the universal cover, Proposition~\ref{prop:universal}. 

These are related as follows: The group $G$ acts freely on the frame bundle of $M$ (by the same argument as in Lemma~\ref{lem:fract}), and $\Fr(M)/G\cong \Fr(\cO)$. 
Associating to an (orbifold) diffeomorphism $f$ of $M$ or $\cO$ the (manifold) diffeomorphism $\Fr(f)$ yields the vertical arrows in the diagram:
\[
\begin{tikzcd}[column sep=scriptsize] \one \arrow[r] & G \arrow[r] \arrow[d] & N_{\Diff(M)}(G) \arrow[r] \arrow[d, "\Fr"] & \Diff(\cO) \arrow[r] \arrow[d, "\Fr"] & \one \\ 
\one \arrow[r] & G \arrow[r] & N_{\Aut(\Fr(M),\theta)}(G) \arrow[r] & \Aut\bigl(\Fr(\cO),\theta\bigr) \arrow[r] & \one \end{tikzcd}
\]
The vertical maps are isomorphisms, which yields the surjectivity of the right horizontal map in the bottom row. We can take the quotient and obtain: 

\begin{theorem}
\label{thm:twoperspectives}
All maps in the diagram below are isomorphisms
\[
\begin{tikzcd}[column sep=large, row sep=large] N_{\Diff(M)}(G)/G \arrow[r] \arrow[d] & \Diff(\cO) \arrow[d] \\ N_{\Aut(\Fr(M),\theta)}\!(G)/G \arrow[r] & \Aut\!\bigl(\Fr(\cO),\theta\bigr). \end{tikzcd}
\]
\end{theorem}

\section{The diffeomorphism group of \texorpdfstring{$\R^n/G$}{Rn/G}}
We now use the previous section to understand the homotopy type of diffeomorphism groups of the local models of orbifolds. For this section we work in the $C^\infty$-compact-open topology. This is different from the Lie-theoretic topology discussed in Section~\ref{sec:lie} and Appendix~\ref{App:Ofd_diff_mfd_struct} if the orbifold is non-compact. 

Let $G$ be a finite group acting linearly and effectively on $\R^n$. By picking an invariant metric, we may assume that $G$ acts by isometries. Thus we will see $G$ as a subgroup of $\mO(n)$. The normalizer $N_{\mO(n)}(G)$ of $G$ in $\mO(n)$ is a finite-dimensional Lie group as it is a closed subgroup of the Lie group $\mO(n)$. 

We will need the following lemma:
\begin{lemma}
\label{lem:homotopyogl}
 Let $G$ be a subgroup of $\mO(n)$. Then the inclusion $N_{\mO(n)}(G)\rightarrow N_{\GL(n)}(G)$ is a $G$-equivariant homotopy equivalence, hence $N_{\mO(n)}(G)/G$ is homotopy equivalent to $N_{\GL(n)}(G)/G$.
\end{lemma}
\begin{proof}
Let $A\in N_{\GL(n)}(G)$ and let $g\in G$. Then $Ag=hA$ for some $h\in G$, and $A^Th=gA^T$ as $g^{-1}=g^T$ and $h^{-1}=h^T$. Thus $A^TAg=gA^TA.$ It follows that $g$ preserves every eigenspace of $A^TA$. As $A^TA$ is positive definite and symmetric, we can define the operator $(A^TA)^s$ by scalar multiplication with $\lambda^s$ on every eigenspace with eigenvalue $\lambda$. It follows that $(A^TA)^sg=g(A^T A)^s$ for every $s\in \R$. In particular $(A^TA)^s\in N_{\GL(n)}(G)$ for all $s$.

Now define the map $H\colon [0,1]\times N_{\GL(n)}(G)\rightarrow N_{\GL(n)}(G)$, by $H(t,A)=A\, (A^T A)^{-\frac{t}{2}}$. Then $H(t,A)\in N_{\GL(n)}(G)$ for all $t$ and $H(1,A)H(1,A)^T=A(A^TA)^{-1}A^T$ is the identity. Finally $H(t,gA)=gH(t,A)$ thus $H$ is an equivariant deformation retraction from $N_{\GL(n)}(G)$ to $N_{\mO(n)}(G)$. 
\end{proof}

\begin{theorem}
\label{thm:alexander}
 The inclusion $N_{\mO(n)}(G)\rightarrow N_{\Diff(\R^n)}(G)$ is a $G$-equivariant homotopy equivalence, hence $\Diff(\R^n/G)$ is homotopy equivalent to $N_{\mO(n)}(G)/G$.
\end{theorem} 

\begin{proof}
For brevity we write $N=N_{\Diff(\R^n)}(G)$ and $N_0=\{F\in N\,|\,F(0)=0\}$, the subgroup of diffeomorphisms that fix the origin. The fixed point set $(\R^n)^G=\{x\in \R^n\,|\,G\, x=x\}$ is a linear subspace hence contractible. Evaluation of a diffeomorphism at $0$ gives the fiber bundle
\[
 N_0\rightarrow N\rightarrow (\R^n)^G
\]
which shows that $N$ is $G$-equivariantly homotopy equivalent to $N_0$. The Alexander trick works equivariantly: Let $H\colon [0,1]\times N_0\rightarrow N_0$ be the map
\[
H(t,F)(x)=\begin{cases}\frac{1}{t}F(tx)\qquad &t>0\\
T_0Fx&t=0.\end{cases}
\]
The map is an equivariant deformation retraction of $N_0$ to $N_{\GL(n)}(G)$. By Lemma~\ref{lem:homotopyogl} this in turn is $G$-equivariantly homotopy equivalent to $N_{\mO(n)}(G)$.
\end{proof}

\subsection{Elementary examples}

The group $N_{\mO(n)}(G)/G$ is a representation-theoretic object and its homotopy type can be understood. Before taking up the general problem we first discuss some elementary examples:

\begin{example}[Rotations in the plane] \label{ex:rot} Every element in $\mO(2)$ is of the form \[R_\theta=\left(\begin{array}{rr}\cos(\theta)&-\sin(\theta)\\ \sin(\theta)&\cos(\theta)\end{array}\right) \qquad \text{or } \qquad S_\theta=\left(\begin{array}{rr}\cos(\theta)&\sin(\theta)\\ \sin(\theta)&-\cos(\theta)\end{array}\right).\] 
Let $\Z_p$ act on the plane via $\Z_p\cong\{R_{2\pi k/p}\,|\,k=0,\ldots, p-1\}$, the subgroup of rotations generated by a rotation by angle $2\pi/p$. We claim that 
\[
N_{\mO(2)}(\Z_p)=\mO(2).
\]
Indeed, every $R_\theta$ obviously commutes with every $R_{2\pi k/p}$, and $S_{\theta}R_{2\pi k/p}=R_{2\pi (p-k)/p}S_{\theta}$. Thus $\Diff(\R^2/\Z_p)\simeq N_{\mO(2)}(\Z_p)/\Z_p=\mO(2)/\Z_p\simeq \mO(2)$.
\end{example}

\begin{example}[Reflection in the plane]
 Let $\Z_2$ now act on the plane via $\Z_2\cong\{S_0,\id\}$, the subgroup generated by the reflection through the $x$-axis (this is a different action from the case $p=2$ of the previous example). As $R_\theta S_0=S_{2\theta}R_{\theta}$, we find that $R_\theta\in N_{\mO(2)}(\Z_2)$ if and only if $\theta=0$ or $\theta=\pi$. A similar argument gives that $S_{\theta}\in N_{\mO(2)}(\Z_2)$ if and only if $\theta=0$ or $\theta=\pi$. Thus 
 \[
 N_{\mO(2)}(\Z_2)=\{R_0,R_\pi,S_0,-S_0\}\cong \Z_2\oplus \Z_2. 
 \]
 It then follows that $\Diff(\R^2/\Z_2)\simeq N_{\mO(2)}(\Z_2)/\Z_2$. Thus $\Diff(\R^2/\Z_2)$ has two contractible components. 
\end{example}

\subsection{Representation-theoretic calculations}

 We compute the homotopy type of the identity component of $N_{\mO(n)}(G)/G$ using the isotypical decomposition of the representation. Following~\cite[Diagram (133)]{bradlowpradagothen} we quotient the exact sequence \eqref{eq:centnorm} by the following exact sequence
 \[
 \mathbf{1}\rightarrow Z(G)\rightarrow G\rightarrow \mathrm{Inn}(G)\rightarrow \mathbf{1}
 \]
 where $Z(G)=\{g\in G\,|\,gh=hg \text{ for all } h\in G\}$ is the \emph{center} of $G$ and $\mathrm{Inn}(G)=\{c_h\in \mathrm{Aut}(G)\,|\,c_h(g)=hgh^{-1} \text{ for } h\in G\}$ is the group of \emph{inner automorphisms} of $G$. For $H=\mO(n)$ this gives
 \[
\mathbf{1}\rightarrow C_{\mO(n)}(G)/Z(G)\rightarrow N_{\mO(n)}(G)/G\rightarrow \mathrm{Out}_{\mO(n)}(G)\rightarrow \mathbf{1}
 \]
 where $\mathrm{Out}
 (G)=\mathrm{Aut}(G)/\mathrm{Inn}(G)$ is the group of \emph{outer automorphisms} and the image of $N_{\mO(n)}(G)/G\rightarrow \mathrm{Out}(G)$ is defined to be $\mathrm{Out}_{\mO(n)}(G)$. This is a finite group. 

 The representation of the finite group $G$ is not assumed to be irreducible. Consider the \emph{isotypical decomposition} of the representation
 \[
\R^n\cong \bigoplus_{i}W_i^{\oplus m_i}
 \]
 
 Here the $W_i$'s denote non-isomorphic irreducible representations, and $m_i$ the multiplicities of these representations. By Schur's lemma each irreducible representation $W_i$ is of real, complex or quaternionic type $D_i=\R,\C,\mathbb H$. We let $\mU(m_i)$ denote the unitary group, and $\mSp(m_i)$ the compact symplectic group. By~\cite[Lemma 3.3 and Lemma 3.4]{druschel} we find that
 \[
 C_{\mO(n)}(G)\cong \prod_{D_i=\R} \mO(m_i)\times\prod_{D_i=\C} \mU(m_i)\times\prod_{D_i=\mathbb{H}} \mSp(m_i)
 \]
The connected component of the identity is
 \[C^0_{\mO(n)}(G)\cong \prod_{D_i=\R} \mSO(m_i)\times\prod_{D_i=\C} \mU(m_i)\times\prod_{D_i=\mathbb{H}} \mSp(m_i).\]
 
\begin{prop}
\label{prop:isotypical}
 Consider the situation above. Then $(N_{\mO(n)}(G)/G)^0$, the connected component of the identity, is isomorphic as a Lie group to
 \[
(N_{\mO(n)}(G)/G)^0 \cong \frac{C^0_{\mO(n)}(G)}{C^0_{\mO(n)}(G)\cap Z(G)}.
 \]
\end{prop} 
\begin{proof}
The group $C_{\mO(n)}(G)/Z(G)$ is both open and closed in $N_{\mO(n)}(G)/G$ as it is the kernel of the continuous homomorphism to the finite group $\mathrm{Out}_{\mO(n)}(G)$. It therefore has the same identity component as $N_{\mO(n)}(G)/G$, which is $\frac{C^0_{\mO(n)}(G)}{C^0_{\mO(n)}(G)\cap Z(G)}$.
\end{proof}

Combining Proposition~\ref{prop:isotypical} with Theorem~\ref{thm:alexander} gives:
\begin{theorem}
\label{thm:diffeochart}
 The connected component of $\Diff(\R^n/G)$ containing the identity is homotopy equivalent to $\frac{C^0_{\mO(n)}(G)}{C^0_{\mO(n)}(G)\cap Z(G)}$.
\end{theorem}

\begin{example}
    The quaternion group $Q=\{\pm 1,\pm i,\pm j,\pm k\}$ acts on $\mathbb{H}\cong \R^4$. This is an irreducible representation of quaternionic type. Thus $C^0_{\mO(4)}(Q)=C_{\mO(4)}(Q)=\mSp(1)$. The center $Z(Q)=\{\pm 1\}\cong \Z_2$. Thus by Theorem~\ref{thm:diffeochart} the identity component of $\Diff(\R^4/Q)$ is homotopy equivalent to $ \mSp(1)/\{\pm 1\}=\mSO(3)$, where we used that the universal cover of $\mSO(3)$ is $\mSp(1)$. 
\end{example}
\section{Homotopy types for two-dimensional orbifolds}

Two-dimensional orbifolds were classified by Thurston~\cite[Chapter 13]{thurston}. An important invariant in the classification is the Euler characteristic of the orbifold. Let $\cO$ be a closed connected two-dimensional orbifold, with $k$ corner reflectors, where the local groups are the dihedral groups of order $2n_1,\ldots,2n_k$, and $l$ elliptic points, where local groups are finite rotation groups of order $m_1,\ldots,m_l$, and any number of mirror reflections. Let $\chi(O)$ be the Euler characteristic of the underlying space. Then the orbifold Euler characteristic is
\[
\chi(\cO)=\chi(O)-\frac{1}{2}\sum_{i=1}^k \left(1-\frac{1}{n_i}\right)-\sum_{i=1}^l\left(1-\frac{1}{m_i}\right)
\]
(Thus the orbifold Euler characteristic does not explicitly depend on the mirror reflections). An orbifold has an \emph{elliptic}, \emph{parabolic} or \emph{hyperbolic} structure depending on whether the Euler characteristic is positive, zero, or negative.  Orbifolds with non-positive Euler characteristic are all good, hence the homotopy type of their diffeomorphism groups can in principle be understood through the isomorphism in Proposition~\ref{prop:universal}. In practice this is of course difficult. Important results in this direction are~\cite{earlemcmullen} and~\cite{earle}. They prove that the connected components of the diffeomorphism groups of closed hyperbolic orbifolds are contractible. 

The orbifolds with positive Euler characteristic are of two types: The bad ones and the good ones (called elliptic by Thurston). We will study the bad orbifolds: They come in two families: Those orbifolds that have underlying spaces the sphere $S^2$ and those that have underlying space the disc $D^2$. The $D^2$ family has a universal cover in the $S^2$ family. The orbifold family with underlying space $S^2$ consists of two types: The teardrops $S^2(p)$ have a single elliptic point of order $p$, and the spindles $S^2(p,q)$ have two elliptic points of orders $p$ and $q$. If $p$ and $q$ are not coprime, the spindle orbifolds are covered by $S^2(p/\gcd(p,q),q/\gcd(p,q))$. We study the diffeomorphism groups of the teardrops and the spindles. 

\subsection{The teardrops}
The goal of this section is to understand the homotopy type of the diffeomorphism groups of the teardrop orbifolds. We will reduce this to a computation of the homotopy type of an embedding space. A codimension-$0$ orbifold embedding is a smooth orbifold map that is a diffeomorphism onto its image. We denote by $\Emb(\cO,\cP)$ the space of orbifold embeddings. 

Let $D^2(p)\subseteq S^2(p)$ be a fixed two-dimensional disc with a $\Z_p$ singularity. A similar argument as in the proof of Theorem~\ref{thm:alexander} gives:
\begin{prop}
\label{prop:embeddings}
 $\Emb(D^2(p),S^2(p))$ is homotopy equivalent to $\mO(2)$.
\end{prop}
\begin{proof}
Every embedding must preserve the singular point. For every embedding $f\in \Emb(D^2(p),S^2(p))$, we let 
\[a(f):=\frac{1}{2}\sup\{t\in[0,1]\,|\,f(tx)\in D^2(p) \text{ for all } x\in D^2(p)\}.\] 
This defines a continuous map $a\colon \Emb(D^2(p),S^2(p))\rightarrow (0,\frac{1}{2}].$ Define the homotopy $H\colon [0,1]\times\Emb(D^2(p),S^2(p))\rightarrow \Emb(D^2(p),S^2(p))$ by $H(t,f)(x)=f((1+t(a(f)-1))x).$ This deforms $\Emb(D^2(p),S^2(p))$ into the subspace $E\subseteq\Emb(D^2(p),S^2(p))$ of embeddings that map $D^2(p)$ into the interior of $D^2(p)$. 

We now pass to a single orbifold chart around the singular point of the form $(\R^2,\Z_p,\varphi)$ containing $D^2(p)$. Thus we view $E$ as a subspace of $\Emb(D^2(p),\R^2/\Z_p)$. The scaling argument shows that these spaces are homotopy equivalent. Every orbifold embedding in $\Emb(D^2(p),\R^2/\Z_p)$ lifts to a map $\w f\colon \w D^2\rightarrow \R^2$ such that $\w f(g \w x)=\alpha(g)\w f(\w x)$ for some automorphism $\alpha\colon \Z_p\rightarrow \Z_p$. The homotopy 
\[
H'(t,\w f)(x)=\begin{cases}\frac{1}{t}(\w f(t x))&t>0\\
T_0 \w fx &t=0
\end{cases}
\]
deformation retracts the space of lifted embeddings to $N_{\GL(2)}(\Z_p)$. Hence $\Emb(D^2(p),S^2(p))$ is homotopy equivalent to $N_{\GL(2)}(\Z_p)/\Z_p$, which by Lemma~\ref{lem:homotopyogl} and Example~\ref{ex:rot} is homotopy equivalent to $\mO(2)$. 
\end{proof}

We can now compute the homotopy type of the diffeomorphism groups of the teardrop orbifolds. 

\begin{prop}
\label{prop:teardrop}
Let $p>1$. Then the orbifold diffeomorphism group $\Diff(S^2(p))$ is homotopy equivalent to $\mO(2)$. 
\end{prop}
\begin{proof}
Restricting the orbifold diffeomorphism to $D^2(p)$ gives the fibration
\begin{equation}
\label{eq:teardrop}
\Diff(S^2(p);D^2(p))\rightarrow \Diff(S^2(p))\rightarrow \Emb(D^2(p),S^2(p)).
\end{equation}
 The space $\Diff(S^2(p);D^2(p))$ is the subgroup of orbifold diffeomorphisms that fix the subspace $D^2(p)$ pointwise. As in the manifold case we find that the space $\Diff(S^2(p);D^2(p))$ is homotopy equivalent to $\Diff (S^2(p)\setminus \mathrm{int}(D^2(p));\partial D^2(p))$ by a collaring argument. The latter group is the (manifold) diffeomorphism group $\Diff(D^2;\partial D^2)$ of diffeomorphisms of the disc that fix the boundary pointwise. The space $\Diff(D^2;\partial D^2)$ is contractible, which is due to Smale~\cite{smale}. Thus $\Diff(S^2(p))$ is homotopy equivalent to $\Emb(D^2(p),S^2(p))$ which by Proposition~\ref{prop:embeddings} is homotopy equivalent to $\mO(2)$.
\end{proof}

\subsection{The spindles}
Let us now investigate the spindles $S^2(p,q)$. 

\begin{prop}
\label{prop:spindle}
 For $p,q>1$ with $p\not=q$ we have that $\Diff(S^2(p,q))$ is homotopy equivalent to $\mO(2)$.
\end{prop}
\begin{proof}
Take two disjoint discs $D^2(p)$ and $D^2(q)$ around the singular points. As $p\not=q$ we can use the scaling trick of Proposition~\ref{prop:embeddings} on each disc at the same time to show that $\Emb(D^2(p)\cup D^2(q),S^2(p,q))$ is homotopy equivalent to $\mO(2)\times \mO(2)$. Note that the restriction map $\Diff(S^2(p,q))\rightarrow \Emb(D^2(p)\cup D^2(q),S^2(p,q))$ is not surjective: An embedding of disjoint discs may reverse the orientation in one disc, and not in the other disc. Embeddings that arise from underlying diffeomorphisms of $S^2(p,q)$ cannot do that. Let $\Emb^c(D^2(p)\cup D^2(q),S^2(p,q))$ denote the coherently oriented embeddings. Then we have a fibration
\begin{equation}
\label{eq:spindle}
\Diff(S^2(p,q);D^2(p)\cup D^2(q))\rightarrow \Diff(S^2(p,q))\rightarrow \Emb^c(D^2(p)\cup D^2(q),S^2(p,q)).
\end{equation}
The fiber is the subgroup of diffeomorphisms fixing both discs pointwise. Again, by a collaring argument we have a homotopy equivalence
\[\Diff(S^2(p,q);D^2(p)\cup D^2(q))\simeq \Diff(S^2(p,q)\setminus \mathrm{int}(D^2(p)\cup D^2(q));\partial D^2(p)\cup \partial D^2(q))
,\]
which is nothing but the (manifold) diffeomorphism group $\Diff(A;\partial A)$ of the annulus fixing the boundary pointwise. This group has contractible components and $\pi_0(\Diff(A;\partial A))\cong\Z$ is generated by the Dehn twists, see for example~\cite{gramain}. The space $\Emb^c(D^2(p)\cup D^2(q),S^2(p,q))
$ can be identified with the space $\{(A,B)\in \mO(2)\times \mO(2)\,|\,\det A=\det B\}$. Thus $\pi_1(\Emb^c(D^2(p)\cup D^2(q),S^2(p,q)))=\Z^2$ and $\pi_0(\Emb^c(D^2(p)\cup D^2(q),S^2(p,q)))=\Z_2$ and all other homotopy groups vanish. The long exact sequence of homotopy groups then shows that $\pi_k(\Diff(S^2(p,q)))=1$ for $k\geq 2$ and 
\[
\mathbf{1}\rightarrow \pi_1(\Diff(S^2(p,q)))\rightarrow \Z^2\xrightarrow{\partial}\Z\rightarrow \pi_0(\Diff(S^2(p,q)))\rightarrow \Z_2\rightarrow \mathbf{1}.
\]
We claim that $\partial$ is surjective. Indeed the generators of $\pi_1(\Emb^c(D^2(p)\cup D^2(q),S^2(p,q)))$ are loops given by a full rotation of one of the discs in the underlying space $S^2$. Lifting this loop to the fiber will produce a Dehn twist on the annulus. 

Thus $\partial$ is surjective and $\pi_1(\Diff(S^2(p,q)))\cong\Z$. The group $\Diff(S^2(p,q))$ has two connected components, which must be homotopy equivalent to each other as it is a topological group. All higher homotopy groups vanish. As $\Diff(S^2(p,q))$ has the homotopy type of a CW complex, it must be homotopy equivalent to $\mO(2)$.
\end{proof}

\begin{remark}
 If $p=q>1$ the orbifold is good. We can follow the arguments in Proposition~\ref{prop:spindle} to compute the homotopy type of the orbifold diffeomorphism group, but there are now two additional connected components of the diffeomorphism group obtained by a reflection of $S^2$ that swaps both singular points. Thus in this case $\Diff(S^2(p,p))\simeq \mO(2)\times \Z_2$.
\end{remark}

\begin{remark}
There is a better description of the homotopy type which takes into account the groups. The spindles and teardrops arise as a quotient of $S^3\subset \C^2$ by the isometric $S^1\subset \C$ action $z(v,w)=(z^pv,z^qw)$ for $q,p$ coprime.  The invariant Riemannian metric on $S^3$ descends to the teardrop and the spindle. The metric is radially symmetric (where the singularities are in the north and south pole). Then it is clear that $\mO(2)$ acts on the teardrops and the spindles. The inclusions $\mO(2)\rightarrow \Diff(S^2(p,q))$ and $\mO(2)\rightarrow \Diff(S^2(p))$ are  homotopy equivalences as they induce isomorphisms on all homotopy groups. 
\end{remark}

\section{Smooth structures on diffeomorphism groups}
\label{sec:lie}
The main result of~\cite{S15} was the construction of a Lie group structure for the group $\Diff (\cO)$ of an effective (possibly non-compact) orbifold. 

The exact sequence \eqref{eq:frameseq} is split by $\Fr$. This means that we can construct $\Aut(\Fr(\cO))$ as a semidirect product of the outer terms $\Diff(\cO)$ and $\Gau(\Fr(\cO))$. We will define a Lie group structure on $\Gau(\Fr(\cO))$ and we use the semidirect product description to give $\Aut(\Fr(\cO))$ a Lie group structure such that all the maps in the sequence are smooth group homomorphisms. 

We will assume that $\cO$ is a compact orbifold. 
The approach is analogous to the constructions for gauge groups of principal bundles, see e.g.~\cite{AaCaMaM89} or~\cite[Section 7.1.1]{GaN26}.

If $M$ is a manifold, then an element $F\in \Gau(\Fr(M))$, 
can be identified with a map $\tau\colon \Fr(M)\rightarrow \GL(n)$. Namely for every $p_x\in \Fr(M)$ we can write $F(p_x)=p_x \tau(p_x)$ for a unique $\tau(p_x)\in \GL(n)$ as $\GL(n)$ acts freely. Equivariance of $F$ then implies that $\tau (p_xA)=A^{-1}\tau(p_x)A$. Conversely, any smooth map $\tau$ satisfying this equation defines an element of $\Gau(\Fr(M))$. In the orbifold case pointwise the element $\tau(p_x)\in \GL(n)$ for $p_x\in \Fr(\cO)$ is not uniquely determined if $\Gamma_x\neq 1$. Indeed, if $p_x\tau(p_x)=p_x\tau'(p_x)$, then $\tau(p_x)(\tau'(p_x))^{-1}\in \GL(n)_{p_x}=\{A\in \GL(n)\,|\, p_x\, A=p_x\}\cong \Gamma_x$. This non-uniqueness is only a pointwise issue: We will show in Proposition~\ref{prop:map2geom_gauge} that we can recover $\tau$ uniquely from $F$ using the solder form. Define the \emph{mapping gauge group} 
\begin{align*}
 \cG= \{\tau\in C^\infty&(\Fr(\cO),\GL(n))\,|\, \tau(p_xA)=A^{-1}\tau(p_x)A,\\
 &\text{ for all } p_x\in \Fr(\cO),\ A\in \GL(n) \}, \end{align*} equipped with its pointwise group structure.

\begin{prop}
\label{prop:map2geom_gauge}
The map
 \begin{align}\label{eq:Psimap} 
\Psi\colon \cG \rightarrow \Gau(\Fr(\cO)), \qquad \Psi(\tau)(p_x):=p_x\tau(p_x) \end{align}
is well defined and a group isomorphism.
\end{prop}
\begin{proof}
We need to check $\Psi(\tau)$ takes values in $\Gau(\Fr(\cO))$. The only non-trivial step is to show that $\Psi(\tau)$ is $\GL(n)$-equivariant. This follows from the computation
 \begin{equation*}
 \Psi(\tau)(p_x A)=(p_xA)\tau(p_xA)=p_xAA^{-1}\tau(p_x)A=\Psi(\tau)(p_x)A.
 \end{equation*}
It is a group homomorphism as
\[
\Psi(\sigma)\circ\Psi(\tau)(p_x)=p_x\tau(p_x)\sigma(p_x\,\tau(p_x))=p_x\tau(p_x)\tau(p_x)^{-1}\sigma(p_x)\tau(p_x)=\Psi(\sigma\tau)(p_x).
\] 
 We can recover $\tau$ from $\Psi(\tau)$ using the solder form. We have that $\Psi(\tau)=\nu(\id,\tau)$ where $\nu\colon \Fr(\cO)\times \GL(n)\rightarrow \Fr(\cO)$ is the $\GL(n)$-action. Then by the product rule
 \[
T_{p_x}\Psi(\tau)=T_{p_x}\nu_{\tau(p_x)}+T_{\tau(p_x)}\nu(p_x,\cdot)\,T_{p_x}\tau.
 \]
 The second term is tangent to the orbit. Thus 
\begin{equation}
\begin{split}
\label{eq:invariance}
(\Psi(\tau)^\ast \theta)_{p_x} &=\theta_{\Psi(\tau)(p_x)}T_{p_x}\Psi(\tau)=\theta_{p_x\tau(p_x)}T_{p_x}\nu_{\tau(p_x)}\\&=(\nu_{\tau(p_x)}^*\theta)_{p_x}=\tau(p_x)^{-1}\theta_{p_x}.
\end{split}
\end{equation}
Here we used Proposition~\ref{prop:solder_orbi} twice: In the second equality we used that $\theta$ vanishes on vectors tangent to the $\GL(n)$-orbits, i.e. part ii) of the proposition and in the last step we used part i). 

Now let $F\in \Gau(\Fr(\cO))$. As $F$ is $\GL(n)$-equivariant and invertible, it maps orbits to orbits diffeomorphically. Therefore $\ker (F^*\theta)_{p_x}=\ker \theta_{p_x}$, by Proposition~\ref{prop:solder_orbi} part ii) as the kernel corresponds to the tangent space to the orbit. As $\theta_{p_x}$ is surjective there is a unique element $\tau_F(p_x)\in \GL(n)$ such that $(F^*\theta)_{p_x}=\left(\tau_F(p_x)\right)^{-1} \theta_{p_x}$. The map $\tau_F\colon \Fr(\cO)\rightarrow \GL(n)$ is clearly smooth. Equivariance of $F$ means that $F\circ \nu_A=\nu_A \circ F$ for all $A\in \GL(n)$ thus
\begin{align*}
(\tau_F(p_xA))^{-1}A^{-1}\theta_{p_x}
&=(\nu_A^*F^*\theta)_{p_x} =(F^*\nu_A^*\theta)_{p_x}=A^{-1}(\tau_F(p_x))^{-1}\theta_{p_x}
\end{align*}
by Proposition~\ref{prop:solder_orbi} again. As $\theta_{p_x}$ is surjective we find that 
$\tau_F(p_xA)=A^{-1}\tau_F(p_x)A$. 
Thus we have a map $\Phi\colon \Gau(\Fr(\cO))\rightarrow \cG$ given by $\Phi(F)=\tau_F$. We claim it is an inverse to $\Psi$. By Equation~\eqref{eq:invariance} we see that $(\Psi(\tau))^*\theta=\tau^{-1}\theta$ hence $\Phi(\Psi(\tau))=\tau$. By Equation~\eqref{eq:invariance}
\[
(\Psi(\Phi(F)))^*\theta=\Phi(F)^{-1}\theta=F^{*}\theta.
\]
Thus $\Psi(\Phi(F))\circ F^{-1}$ preserves $\theta$. By Proposition~\ref{prop:orb2} $\Psi(\Phi(F))\circ F^{-1}=\Fr(\id_\cO)=\id_{\Fr(\cO)}$. Thus $\Psi(\Phi(F))=F$, and $\Psi$ and $\Phi$ are each other's inverses. 
\end{proof}

To obtain a Lie group structure on $\Gau(\Fr(\cO))$ we construct a Lie group structure on $\cG$. The following construction is standard~\cite[Section 10]{GaN26}, but as $\Fr(\cO)$ is non-compact, even for compact $\cO$, we have to be careful.

Let $E=\text{Mat}_n(\R)$ be the matrix algebra of $n\times n$-matrices. Since $\Fr(\cO)$ is $\sigma$-compact and $E$ is a Fr\'{e}chet space the smooth mapping space $C^\infty (\Fr(\cO),E)$ is also a Fr\'{e}chet space with respect to the pointwise operations and the compact-open $C^\infty$-topology (cf.~e.g.~\cite[Section 3]{GaS22}). 
Now $E$ is a Banach algebra, and the pointwise multiplication of elements turns $C^\infty (\Fr(\cO),E)$ into a locally convex algebra by~\cite[Lemma 12.2.2 (a)]{GaN26}.
While $\GL(n)\subseteq E$ is open, as $\Fr(\cO)$ is non-compact, the unit group $C^\infty (\Fr(\cO),\GL(n))$ of $C^\infty (\Fr(\cO),E)$ is not open in the compact-open $C^\infty$-topology. Indeed, since this topology only controls values on compact sets, every open neighborhood of a map with values only in the open set $\GL(n)$ contains maps taking some values in $E\setminus \GL(n)$. Therefore $C^\infty (\Fr(\cO),\GL(n))$ does not directly inherit an infinite-dimensional Lie group structure from $C^\infty (\Fr(\cO),E)$.

Equivariance saves us however. The conjugation $c_A \colon E \rightarrow E, c_{A}(e)=A^{-1}eA$ for $A\in \GL(n)$ is continuous. Since point evaluations are continuous linear in the compact-open $C^\infty$-topology, see~\cite[Example 4.1(a) and Lemma 4.3 (a)]{GaS22}, we obtain a closed linear subspace 
\begin{align}\label{eq:ad-inv_subalg}
\begin{split}\mathcal{A}:= \{\tau \in C^\infty (\Fr(\cO),E)\,|&\,\tau(p_xA)=A^{-1}\tau(p_x)A, \\
&\text{ for all } p_x\in \Fr(\cO), A\in \GL(n)\}
\end{split}
\end{align}
which is a (Fr\'{e}chet) subalgebra of $C^\infty (\Fr(\cO),E)$, so in particular a locally convex algebra. 
Moreover, we have the following:
\begin{lemma}\label{lem:properties_A}
 The algebra $\mathcal{A}$ is a continuous inverse algebra which is locally $m$-convex, i.e. its topology is generated by submultiplicative seminorms $q(xy)\leq q(x)q(y),$ for all $x,y \in \mathcal{A}$.
\end{lemma}
\begin{proof}
The unit group of $\mathcal{A}$ is $\cG$. This is an open subset: As $\cO$ is compact, the $\GL(n)$-action on $\Fr(\cO)$ is cocompact. Thus there exists a compact set $C\subset \Fr(\cO)$ such that $C\GL(n)=\Fr(\cO)$, which one can construct using local charts. By equivariance it follows that if $\tau(C)\subset \GL(n)$ for $\tau\in \mathcal A$, then the full image of $\tau$ is contained in $\GL(n)$. Thus
\[
\cG=\{\tau\in \mathcal{A}\,|\,\tau(C)\subseteq \GL(n)\},
\]
which is open in the compact-open topology. Thus $\mathcal{A}$ is a locally convex algebra whose unit group is open and inversion is continuous (as it is given by postcomposition with inversion in the matrix algebra) on this set. Thus $\mathcal{A}$ is a continuous inverse algebra.

For the local $m$-convexity we were not able to find a proof readily citable in the literature. The property follows from the Leibniz rule. We give the argument for the reader's convenience. Choose a submultiplicative norm $\lVert \cdot\rVert_E$ on $E$.
On a compact subset $K\subseteq \Fr(\cO)$ contained in a coordinate chart $(U,\varphi)$ (which we suppress in the notation in the following), define 
\[ q_{K,k}(f) := \sum_{|\alpha|\leq k} \frac{1}{\alpha!} \sup_{x\in K} \lVert \partial^\alpha f(x)\rVert_E, \]
where we use standard multiindex notation, e.g. for $\alpha=(\alpha_1,\ldots, \alpha_d)$
\[ \alpha!:=\alpha_1!\cdots\alpha_d!, \quad \binom{\alpha}{\beta} := \prod_{j=1}^d\binom{\alpha_j}{\beta_j}, \quad \binom{\alpha_j}{\beta_j} = \frac{\alpha_j!}{\beta_j!(\alpha_j-\beta_j)!}. \]
Up to the weights $1/\alpha!$, these are the standard seminorms for the compact-open $C^\infty$-topology.  Thus the seminorms $q_{K,k}$ with $K\subseteq \Fr(\cO)$ running through all compact subsets contained in some chart domain and $k\in\N_0$ generate the compact-open $C^\infty$-topology on $C^\infty (\Fr(\cO),E)$.
The multi-index Leibniz rule for the multiplication in $E$ (use e.g. inductively~\cite[Exercises 1.2.3 and 1.4.2]{GaN26}) gives 
$\partial^\alpha(fg) = \sum_{\beta\leq\alpha} \binom{\alpha}{\beta} (\partial^\beta f) (\partial^{\alpha-\beta}g)$.
By submultiplicativity of the norm and the identity $\frac{1}{\alpha!}\binom{\alpha}{\beta} = \frac{1}{\beta!(\alpha-\beta)!}$, we obtain
\begin{align*} 
q_{K,k}(fg) &\leq \sum_{|\alpha|\leq k} \sum_{\beta\leq\alpha} \frac{ \sup_K\lVert\partial^\beta f\rVert_E \sup_K\lVert\partial^{\alpha-\beta}g\rVert_E }{ \beta!(\alpha-\beta)! } \\ &\leq \left( \sum_{|\beta|\leq k} \frac{\sup_K\lVert\partial^\beta f\rVert_E}{\beta!} \right) \left( \sum_{|\gamma|\leq k} \frac{\sup_K\lVert\partial^\gamma g\rVert_E}{\gamma!} \right) = q_{K,k}(f)\,q_{K,k}(g). 
\end{align*}
Consequently, $C^\infty(\Fr(\cO),E)$ is locally $m$-convex. 
The seminorms of $C^\infty(\Fr(\cO),E)$ restrict to seminorms generating the subtopology of $\mathcal{A}$, and therefore also $\mathcal{A}$ is locally $m$-convex.
\end{proof}

Before we continue, recall first (e.g. from~\cite{S23,GaN26}) that a Lie group $G$ is called \emph{semiregular} (in the sense of Milnor) if for every smooth curve $\eta \colon [0,1] \rightarrow \Lf(G)$ into its Lie algebra, the Lie type differential equation 
\[\dot{\gamma} (t)=\gamma(t).\eta(t):= T\lambda_{\gamma(t)}(\eta(t)) \qquad \gamma(0)=\one_G,\]
where $\lambda_g(h):=gh$ is left multiplication, has a unique solution $\gamma_\eta \in C^\infty ([0,1],G)$. If moreover, the map $\evol \colon C^\infty([0,1],\Lf(G)) \rightarrow G, \eta \mapsto \gamma_\eta(1)$ is smooth, we say that $G$ is \emph{regular} (in the sense of Milnor).

\begin{theorem}\label{thm:map_gp_Lie}
The unit group $\cG$ of $\mathcal{A}$ is a regular Fr\'{e}chet Lie group. Moreover, $\cG$ is even analytic and a BCH-Lie group, so in particular it is a locally exponential Lie group whose Lie group exponential is given by the push-forward of the Lie group exponential $\exp_{\GL(n)}$.
\end{theorem}
\begin{proof}
The unit group of a continuous inverse algebra is an analytic BCH-Lie group with the claimed exponential map by~\cite[Theorem 10.1.15]{GaN26}. As Lemma~\ref{lem:properties_A} establishes that $\mathcal{A}$ is locally $m$-convex, regularity in the sense of Milnor follows by~\cite[Theorem 10.6.1 and Remark~10.6.2(b)]{GaN26}. 
\end{proof}

Let $f\in \Diff(\cO)$ and $\tau\in \cG$. We compute that
\[
\tau(\Fr(f^{-1})(p_xA))=\tau(\Fr(f^{-1})(p_x)A)=A^{-1}\tau( \Fr(f^{-1})(p_x))A
\]
as $\Fr(f^{-1})$ is equivariant. Thus $\tau\circ \Fr(f^{-1})\in \cG$. This formula therefore defines a left group action $\beta\colon \Diff(\cO)\times \cG\rightarrow \cG$ via $\beta(f,\tau)=\tau \circ\Fr(f^{-1}).$ The semidirect product $\cG\rtimes\Diff(\cO)$ is the group with underlying set $\cG\times \Diff(\cO)$ with the product $(\tau,f)(\sigma,g)=(\tau\,(\beta(f,\sigma)),f\circ g)$. 
\begin{prop}
 The action $\beta$ is a smooth Lie group action, and $\cG\rtimes_\beta\Diff(\cO)$ is naturally a regular Fr\'echet Lie group. 
\end{prop}
\begin{proof}
Similar to the left action $\beta$ we have a right action 
\[\alpha \colon \cG \times \Diff(\cO) \rightarrow \cG,\quad \alpha(\tau , f)=\beta(f^{-1},\tau)=\tau\circ\Fr(f).\]
Note that the right action $\alpha$ induces the left action $\beta$ by composing it with the inversion in $\Diff(\cO)$ (which is smooth as inversion in a Lie group).

The group $\cG$ is the unit group of a continuous inverse algebra $\mathcal{A}$ which is a closed subspace of the Fr\'{e}chet space $C^\infty (\Fr(\cO),E)$ with the compact-open $C^\infty$-topology. By~\cite[Lemma 1.25]{S23}, the map $\alpha$ will be smooth if and only if it is smooth as a mapping to $C^\infty (\Fr(\cO),E)$. Hence~\cite[Lemma 3.13]{GaS22} shows that smoothness of $\alpha$ is equivalent to the smoothness of the adjoint map 
\begin{align}\label{alpha_adj}
\alpha^\wedge \colon \cG \times \Diff(\cO) \times \Fr(\cO) \rightarrow E,\quad \alpha^\wedge(\tau, f,p_x):= \tau(\Fr(f)(p_x))
\end{align}
By~\cite[Lemma 4.3 (a)]{GaS22} $\ev \colon C^\infty (\Fr(\cO),E) \times \Fr(\cO) \rightarrow E, \ev(f,p_x)=f(p_x)$ is smooth. The smooth inclusion $\text{inc} \colon \cG \rightarrow C^\infty (\Fr(\cO),E)$ allows us to rewrite $\alpha^\wedge$ using the smooth map $\Xi$ from Lemma~\ref{smooth_adjoint_section} as 
$\alpha^\wedge (\tau,f,p_x) = \ev(\text{inc}(\tau),\Xi(f,p_x))$. Thus the map~\eqref{alpha_adj} and as a consequence $\alpha$ and $\beta$ are smooth. Then the Lie group structure on the semidirect product can be constructed as in finite-dimensional Lie theory by the usual arguments, cf.~\cite[Section 5.2.2]{GaN26}. Moreover both $\cG$ and $\Diff(\cO)$ are regular Fr\'{e}chet Lie groups (by Theorem~\ref{thm:map_gp_Lie} and~\cite[Theorem 5.35]{S15}). Since regularity is an extension property,~\cite[Theorem 38.6]{KaM97}, the semidirect product is a regular Fr\'{e}chet Lie group.
\end{proof}

We equip $\Gau(\Fr(\cO))$ with the Fr\'echet Lie group structure by declaring $\Psi$ to be an isomorphism of Fr\'echet Lie groups.  We compute that
\begin{align*}
\Psi(\beta(f,\tau))(p_x)&=p_x \tau(\Fr(f^{-1})p_x)=\Fr(f)\left(\Fr(f^{-1})(p_x) \tau(\Fr(f^{-1})p_x)\right)\\
&=\Fr(f)\circ\Psi(\tau)\circ \Fr(f^{-1})(p_x).
\end{align*} 
Thus $\beta'\colon \Diff(\cO)\times \Gau(\Fr(\cO))\rightarrow \Gau(\Fr(\cO))$, $\beta'(f,F)=\Fr(f)\circ F\circ \Fr(f^{-1})$ is a smooth group action as it intertwines the group action $\beta$ under the isomorphism $\cG\rightarrow \Gau(\Fr(\cO))$. Thus $\Gau(\Fr(\cO))\rtimes_{\beta'} \Diff(\cO)$ is Lie group isomorphic to $\cG\rtimes_\beta \Diff(\cO)$. 

The map
$\gamma\colon \Gau(\Fr(\cO))\rtimes_{\beta'}\Diff(\cO)\rightarrow \Aut(\Fr(\cO))$ defined by $\gamma(F,f)=F\circ\Fr(f)$ is a group isomorphism. We equip $\Aut(\Fr(\cO))$ with the regular Fr\'{e}chet Lie group structure by declaring $\gamma$ to be an isomorphism of Lie groups. We then have our main theorem.
\begin{theorem}
 The maps in the exact sequence
 \[
 \one\to\Gau(\Fr(\cO))\to\Aut(\Fr(\cO))\to\Diff(\cO)\to\one
 \]
 are Lie group homomorphisms and the sequence is split by the Lie group homomorphism $\Fr$. Thus $\Aut(\Fr(\cO))$ is a trivial principal $\Gau(\Fr(\cO))$ bundle over $\Diff(\cO)$.
\end{theorem}

The Lie group structure on $\Aut(\Fr(\cO))$ is well known in the case that the compact orbifold $\cO$ is a compact manifold $M$, i.e. there are no singularities. Then the Lie group structure constructed on $\Diff(\cO)$ in~\cite{S15} coincides with the standard Lie group structure of the diffeomorphism group $\Diff(M)$. In this case we get a principal bundle $\Fr(M)\rightarrow M$ and the construction reduces to the structure discussed in~\cite{AaCaMaM89}. Namely we have the following statement:

\begin{remark}\label{rem:old_bundle}
Let $M$ be a compact smooth manifold.
Then cf.~\cite[p. 218 and Corollary 4.4]{AaCaMaM89}, we obtain a short exact sequence 
\[\mathbf{1} \rightarrow \Gau(\Fr(M))\rightarrow \mathrm{Aut} (\Fr(M)) \rightarrow \Diff(M) \rightarrow \mathbf{1}\]
Frame prolongation defines a global section 
$\Fr\colon \Diff(M)\rightarrow\mathrm{Aut}(\Fr(M))$ and $ \cG\rtimes_\beta\Diff(M) \rightarrow \mathrm{Aut}(\Fr(M)), (\tau,f) \mapsto \Psi(\tau)\circ\Fr(f)$, where $\beta(f,\tau) = \tau\circ\Fr(f)^{-1}$, 
is a Lie group isomorphism if we endow the right-hand side with the Lie group structure constructed in~\cite{AaCaMaM89}.
\end{remark}

The construction of the Lie group structure on the automorphism group used that the underlying orbifold is compact. A generalised construction (requiring a more involved topological treatment of the gauge group) should be possible for non-compact orbifolds. We expect that arguments such as those developed in~\cite{Sch13} for the gauge group of a principal bundle over a non-compact manifold can be adapted to the orbifold frame-bundle. This would lead to a Lie group theoretic splitting based on (LF)-spaces instead of the Fr\'{e}chet spaces considered here.

\appendix

\section{The Lie group of orbifold diffeomorphisms}\label{App:Ofd_diff_mfd_struct}

In several arguments in the main text information on charts and the manifold structure of the orbifold diffeomorphism group $\Diff(\cO)$ is needed. We recall the necessary information from~\cite{S15} in this appendix.
As in the rest of the article $\cO$ will denote an effective orbifold, which is paracompact as a topological space but need not be compact.

Let us first recall that the orbifold diffeomorphisms are modelled on the space of orbisections which is topologised as follows

\begin{setup}
The space of orbisections is denoted by $\mathfrak{X}_{\text{Orb}}(\cO)$ and we write $\mathfrak{X}_{\text{Orb},c}(\cO)$ for the space of compactly supported orbisections. If $(\widetilde{U},\Gamma,\psi)$ is an orbifold chart, then every orbisection $\sigma$ admits a canonical lift $\sigma_{\widetilde{U}}\in \mathfrak{X}(\widetilde{U})$ and these lifts are related via the change of charts morphisms to each other.  Thus the space of all compactly supported orbisections can be topologised by choosing a locally finite orbifold atlas $\mathcal{A}:=\{(\widetilde{U}_i,\Gamma_i,\psi_i)\}_{i\in I}$ and declaring the mapping
\begin{align}\label{eq:Lambda_map}
\Lambda_{\mathcal{A}} \colon \mathfrak{X}_{\text{Orb},c} (\cO) \rightarrow \bigoplus_{i\in I} \mathfrak{X}(\widetilde{U}_i), \quad \sigma \mapsto (\sigma_i)_{i\in I}
\end{align}
a (linear) topological embedding, where the spaces of vector fields are endowed with the compact-open $C^\infty$-topology. Note that this map takes its image in the direct sum as the atlas is locally finite.
\end{setup}
One can show,~\cite[Section 3.3]{S15} that the topology on the orbisections does not depend on any choice for the locally finite atlas $\mathcal{A}$. Further, the image of $\Lambda_{\mathcal{A}}$ is a closed subspace of all vector fields related by the orbifold change of charts. 

Recall from~\cite[Section 4]{S15} that every paracompact orbifold admits a Riemannian orbifold metric with associated (Riemannian) orbifold exponential map
$
\exp_{\text{Orb}} \colon T\cO \supseteq \Omega \rightarrow \cO
$
The lifts of $\exp_{\text{Orb}}$ are given by the exponential map of equivariant Riemannian metrics lifting the orbifold Riemannian metric. Moreover, the lifts are related to each other via the orbifold change of charts. In the following we will always fix a Riemannian orbifold metric (all constructions are independent of the actual choice, cf.~\cite[Lemma 5.14]{S15}).

The orbisections form the modelling space of the orbifold diffeomorphism group $\Diff(\cO)$. We recall from~\cite[Section 5.1]{S15} the relevant data used to construct a manifold chart around the identity in $\Diff(\cO)$.

\begin{setup}\label{setup:painful_details}
First an open set $\mathcal{M}\subseteq \mathfrak{X}_{\text{Orb},c}(\cO)$ is constructed from the following data: 
\begin{itemize}
\item there are two locally finite orbifold atlases $\mathcal{A}=\{(\widetilde{U}_i,\Gamma_i,\psi_i)\}_{i\in I},\mathcal{B}:=\{(\widetilde{W}_j,\Gamma_j,\varphi_j)\}_{j\in J}$ and a mapping $\alpha \colon I\rightarrow J$, such that $\overline{U_i}\subseteq W_{\alpha(i)}$.
\item We may assume that every $W_j, j\in J$ has compact closure in $\cO$, so only finitely many $U_i$ meet the closure. 
\item For every $i \in I$ there is a compact subset $\varphi^{-1}_{\alpha(i)}(\overline{U_i})\subseteq K_i^\circ \subseteq K_i \subseteq \widetilde{W}_{\alpha (i)}$.
\item For every $j\in J$ we have a smooth lift $\exp_j \colon T\widetilde{W}_j \supseteq \Omega_j \rightarrow \widetilde{W}_j$ of the orbifold exponential.
\item Shrinking $\Omega_j$ if necessary, we may assume that $(\pi_{\widetilde{W}_j},\exp_j)\colon \Omega_j \rightarrow \widetilde{W}_j \times \widetilde{W}_j$ is an embedding whose open image contains the diagonal.
\end{itemize}
Now $\mathcal{M}$ is constructed with the following properties:
For $\sigma \in \mathcal{M}$ the canonical lift $\sigma_{\alpha(i)} \in \mathfrak{X}(\widetilde{W}_{\alpha(i)})$ satisfies $\sigma_{\alpha(i)}(K_i) \subseteq \Omega_{\alpha(i)}$. 

Therefore one can construct a chart of $\Diff(\cO)$, cf.~\cite[Proposition 5.4]{S15}
\begin{align}\label{mfd_chart:OFGP}
E \colon \mathcal{M} \rightarrow \Diff(\cO), \sigma \mapsto \exp_{\text{Orb}} \circ \sigma
\end{align}
such that for the diffeomorphism $E(\sigma)$ the canonical lift in a chart $(\widetilde{U}_i,\Gamma_i,\psi_i) \in \mathcal{A}$, is (up to an embedding) given by $\exp_{\alpha(i)} \circ \sigma_i$. In particular $E(\sigma)(\overline{U_i})\subseteq W_{\alpha (i)}$ for every $\sigma \in \mathcal{M}$ and $i\in I$.
\end{setup}
We recall from~\cite[Lemma 5.13]{S15} that the manifold structure does not depend on the concrete selection of $\mathcal{A},\mathcal{B}$ nor the concrete orbifold Riemannian metric as long as they satisfy the requirements of~\cite[5.1]{S15}.

Using the notation for charts and compact subsets introduced above, write 
\[\lfloor K,\Omega\rfloor :=\{f \in \mathfrak{X}(\widetilde{W})\,|\,f(K)\subseteq \Omega\}\]
(or obvious variants of this for other compact and open sets).  Note that by convention $e^\sigma := \exp \circ \sigma|_K$ makes sense for all $\sigma \in \lfloor K,\Omega\rfloor $ (where $\Omega$ is the domain of the exponential map).

\begin{lemma}\label{smooth_adjoint_section}
Let $\cO$ be an effective orbifold. Then the following map is smooth 
\begin{align}\label{map_diff_lifting_smooth}
\Xi \colon \Diff(\cO) \times \Fr(\cO) \rightarrow \Fr(\cO), \quad (f,p_x) \mapsto \Fr(f)(p_x).
\end{align} 
\end{lemma}

\begin{proof}
It suffices to show that there is an open neighborhood $\mathcal{U} \subseteq \Diff(\cO)$ of the identity diffeomorphism, such that $\Xi|_{\mathcal{U} \times \Fr(\cO)}$ is smooth. This follows from the fact that $\Fr$ is a group morphism: If $f_0 \in \Diff(\cO)$, then for $ h\in f_0\mathcal{U}$ we have $\Xi(h,p_x)=\Xi(f_0f_0^{-1}h,p_x)=\Fr(f_0f_0^{-1}h)(p_x)=\Fr(f_0) (\Xi (f_0^{-1}h,p_x))$. Now $\Fr(f_0)$ is a diffeomorphism of $\Fr(\cO)$ and $f_0^{-1}h \in \mathcal{U}$.  Thus the map is smooth on the open $f_0$-neighborhood $f_0\mathcal{U}$.

We work with the manifold chart $E \colon \mathcal{M}\rightarrow E(\mathcal{M})\subseteq \Diff(\cO)$ from \eqref{mfd_chart:OFGP}. 
To shrink $\mathcal{M}$ suitably, let $\mathcal{B} =\{(\widetilde{W}_j,\Gamma_j,\varphi_j)\}_{j\in J}$ be the orbifold atlas with canonical map $\Lambda_{\mathcal{B}}$ mapping orbisections to the family of their lifts with respect to $\mathcal{B}$, cf.~\eqref{eq:Lambda_map}. The map $\Lambda_\mathcal{B}$ is continuous linear (where the spaces in the direct sum are endowed with the compact-open topology, cf.~\cite[Section 3.3]{S15}), thus smooth. Its image consists of families of vector fields which are compatible with respect to orbifold change of charts.

We recall from~\ref{setup:painful_details} that for every $i\in I$ there is a compact set $K_i \subseteq \widetilde{W}_{\alpha(i)}$ such that the images of the interiors of the $K_i$ cover $\cO$ and the lift of $E(\sigma)$ is on $K_i$ given by $\exp_{\alpha(i)} \circ\sigma_{\alpha (i)}$. Note that $\Lambda_{\mathcal{B}}$ sends an orbisection in $\mathcal{M}$ to the local representatives on charts in $\mathcal{B}$, i.e. vector fields such that $\sigma_{\alpha(i)}(K_i)\subseteq \Omega_{\alpha(i)}$.

Consider now $p=(x,V)\in \Fr(\cO)$. Then there is (at least) one $i \in I$ such that $x\in \varphi_{\alpha(i)}(K_i^\circ)$. We pick a preimage $x_{\alpha(i)} \in K_i^\circ \cap \varphi^{-1}_{\alpha(i)}(x)$ (the concrete choice is immaterial due to equivariance of the mappings). We write shorter $e^{\sigma}:=\exp_j \circ \sigma$. Choose an open $0$-neighborhood 
\begin{align}
A_i \subseteq \lfloor K_i,\Omega_{\alpha (i)}\rfloor \cap \mathrm{pr}_{\alpha(i)} \circ \Lambda_{\mathcal{B}} (\mathcal{M}) \subseteq \mathfrak{X}(\widetilde{W}_{\alpha(i)}) \text{ and define } \notag \\ H_i \colon A_i \times K_i^\circ \rightarrow \widetilde{W}_{\alpha (i)}, H_i (\sigma, y) := \exp_{\alpha (i)}(y,\sigma (y))=\exp_{\alpha (i)}(y,\mathrm{ev}(\sigma,y)). \label{H_i_map}
\end{align} 
Then \eqref{H_i_map} makes sense and is smooth as the joint evaluation map is smooth by~\cite[Lemma 4.3 (a) and Example 4.1]{GaS22}. 
Further, \eqref{H_i_map} is equivariant in $y$ and as the chart $E$ takes its values in the orbifold diffeomorphisms, if $\sigma_{\alpha(i)} = \text{pr}_{\alpha(i)} \circ \Lambda_\mathcal{B}(\sigma)$ for some $\sigma \in \mathcal{M}$, then the lift $\exp_{\alpha(i)}\circ \mathrm{pr}_{\alpha(i)}\circ \Lambda_{\mathcal{B}}(\sigma) $ of $E(\sigma)$ is a local diffeomorphism. In particular, the derivative $D_yH_i(\sigma, y)$ is a linear isomorphism and thus transforms a frame to a frame.  By the chain rule, the map $(\sigma, y)\mapsto D_y H_i(\sigma,y)$ is smooth in $y$ and $\sigma$ if the map sending $\sigma$ to its derivative is smooth. Working in a local trivialisation of the tangent bundles over an open set $V$, we see that the derivative map is smooth in $\sigma$ as it corresponds to the map $D\colon C^\infty (V,\R^n) \rightarrow C^\infty (V,L(\R^n,\R^n)), F \mapsto  DF $ is smooth, which is true by~\cite[42.19]{KaM97}. 

Hence for every $(\sigma, (y,Z)) \in (\text{pr}_{\alpha(i)} \circ \Lambda_{\mathcal{B}})^{-1}(O_{i}) \times \Fr(\widetilde{U}_i)$ we can rewrite $\Xi(E(\sigma), \Fr(\varphi_{\alpha(i)})(y,Z))$ as the composition
\begin{align}\label{eq:loc_form_Xi}
 \Fr(\varphi_{\alpha (i)})(H_i(\text{pr}_{\alpha(i)} (\Lambda_\mathcal{B}(\sigma)),y), D_y H_i(\text{pr}_{\alpha(i)} (\Lambda_\mathcal{B}(\sigma)),y)Z), \end{align}
where $\Fr(\varphi_{\alpha(i)})\colon \Fr(\widetilde{W}_{\alpha(i)}) \rightarrow \Fr(W_{\alpha(i)})$ is the covering map which is a submersion and $\Fr(\varphi_{\alpha(i)})(\Fr(\widetilde{U}_i))=\Fr(\cO)|_{U_i}$. Note that \eqref{eq:loc_form_Xi} is given as a composition of smooth (in $(\sigma,y)$) maps. Further the composition is given by the equivariant (in $y$) map \eqref{H_i_map} and the map $D_y H_i$ which is equivariant in $(y,Z)$ with respect to the derived action of the local group by construction. For other choices of $x_{\alpha(i)}$ and the chart, these choices are related by orbifold change of charts. Hence \eqref{eq:loc_form_Xi} is independent of the choices of $x_j$ and descends to a smooth map on $\Fr(\cO)|_{U_i}$. By the universal property of the submersion $S_i:=(E,\Fr(\varphi_{\alpha (i)}))\colon \mathcal{M} \times \Fr(\w{U}_i) \rightarrow \Diff(\cO) \times \Fr(\cO)|_{U_i}$, $\Xi$ is smooth on the open $(\id_{\cO},x)$-neighborhood $S_i\left(\mathcal{M} \cap\left(\text{pr}_{\alpha(i)}\circ \Lambda_{\mathcal{B}}\right)^{-1}(A_i) \times \Fr(U_i)\right)$. Note that we can repeat the construction for every $i\in I$. Then we define 
\[\mathcal{V}_j:=\begin{cases} \mathfrak{X}(\tilde{W}_j), & \alpha^{-1}(j)=\emptyset \\
\bigcap_{i \in \alpha^{-1}(j)}A_i, & \alpha^{-1}(j)\neq \emptyset
\end{cases} \qquad \mathcal{V}:=\bigoplus_{j \in J} \mathcal{V}_j.\]
Note that for every $j$ the set $\mathcal{V}_j \subseteq \mathfrak{X}(\tilde{W}_j)$ is open as $\alpha^{-1}(j)$ is always a finite set. Now as the domains $U_i$ cover $\cO$, it suffices to observe that the box neighborhood $\mathcal{V}$ in the locally convex sum is an open neighborhood of the zero-section. Set $\mathcal{U}=E(\mathcal{M}\cap\Lambda_{\mathcal{B}}^{-1}(O))$ and note that this set is open in $\Diff(\cO)$ and contains in particular $\id_{\cO}$. By construction $\Xi$ is smooth on the open set $\mathcal{U} \times \Fr(\cO)$, thus it is smooth.
\end{proof}

\subsection{Relation to homeomorphism groups}
Recall first that for a topological space $X$, the compact-open topology on $C(X,X)$ is generated by the subbase consisting of the neighborhoods
\[\lfloor K,U\rfloor :=\{f\in C(X,X):f(K)\subseteq U\} , K\subseteq X \text{ compact and }U\subseteq X \text{ open,}\]
The following is well known (cf.~\cite{Dij05}): For a Hausdorff topological space $X$ which is locally compact and locally connected, the compact-open topology turns the group of homeomorphisms $\mathrm{Homeo}(X)$ into a topological group.

Thus in particular, if $X$ is a finite-dimensional manifold or orbifold, the compact-open topology turns the homeomorphism group into a topological group.

The following statement is what one would naturally expect, but was not proved in~\cite{S15}.
\begin{lemma}
The natural inclusion $\Diff (\cO) \rightarrow \mathrm{Homeo} (\cO)$ is a morphism of topological groups, i.e. the manifold topology is finer than the (restriction of) the compact-open topology.
\end{lemma}
\begin{proof}
Both $\Diff(\cO)$ and $\mathrm{Homeo}(\cO)$ are topological groups in their respective topologies.
As the canonical inclusion sends an orbifold diffeomorphism to its underlying homeomorphism (see e.g.~\cite[Corollary 2.13]{S15}), it is a group morphism and thus it suffices to prove that it is continuous at the identity element, i.e. that the subbasic neighborhoods of the identity of the form 
$\lfloor K,U\rfloor$, for $K\subseteq \cO$ compact and $U \subseteq \cO$ open, are open in the manifold topology. For any such neighborhood note that $K\subseteq U \subseteq \cO$.

As the construction of the manifold structure does not depend on the choice of orbifold atlases, we can choose new atlases $\widetilde{\mathcal{A}}=\{(\widetilde{U}_i,\Gamma_i,\psi_i)\}_{i\in I},\widetilde{\mathcal{B}}=\{\widetilde{W}_j,\Gamma_j,\varphi_j)\}_{j \in J}$ such that $U_i\cap K \neq \emptyset$ implies $W_{\alpha(i)} \subseteq U$.
To see this pick for each $x \in K$ an orbifold chart domain $x \in U_x \subseteq W_x \subseteq U$. Finitely many of these cover $K$ by compactness. We expand the two finite families  to a locally finite atlas of $\cO$. Then shrink the zero-neighborhood of orbisections such that for each of the finitely many $i$ of the $U_i$ covering $K$ we have $E(\sigma)(U_i)\subseteq W_{\alpha (i)} \subseteq U$. 
This yields an open neighborhood $\widetilde{\mathcal{M}} \subseteq \mathfrak{X}_{\text{Orb},c}(\cO)$ of the zero-orbisection such that $E(\widetilde{\mathcal{M}} \cap \mathcal{M}) \subseteq \lfloor K,U\rfloor$. Thus every compact-open neighborhood of the identity contains an identity neighborhood for the manifold topology. Hence the inclusion is continuous. 
\end{proof}

\begin{prop}\label{prop:ofd_emb_fO}
The map $F \colon \Diff(\cO) \rightarrow \mathrm{Homeo}(\Fr(\cO)), \phi \mapsto \Fr(\phi)$ is a morphism of topological groups.
\end{prop}

\begin{proof}
Note first that \Cref{thm:Fr_gp_mor} shows that $F$ is a group morphism (as the diffeomorphism group is a subgroup of $\mathrm{Homeo}(\Fr(\cO))$. Computing the adjoint map to $F$, we see that $F^\wedge \colon \Diff(\cO) \times \Fr(\cO) \rightarrow \Fr(\cO)$ coincides exactly with the mapping $\Xi(\phi ,p)=\Fr(\phi)(p)$ which is smooth by Lemma~\ref{smooth_adjoint_section}. Since $\mathrm{Homeo}(\Fr(\cO)) \subseteq C(\Fr(\cO),\Fr(\cO))$ carries the subspace topology we can invoke the exponential law for the compact open topology to establish continuity.  As the adjoint map $F^\wedge$ is thus continuous, we deduce from the fact that $\Fr(\cO)$ is locally compact and the exponential law for the compact-open topology~\cite[Proposition B.2.8]{S23}, that $F$ is continuous.
\end{proof}

\begin{remark}
As $\Fr$ takes its image inside of the subgroup $\Diff (\Fr(\cO)$ of smooth diffeomorphisms in the homeomorphisms, an argument similar to the proof of \Cref{prop:ofd_emb_fO} shows that the mappings $F_k \colon \Diff(\cO) \rightarrow \Diff^{C^k}(\Fr(\cO)), \Phi \mapsto \Fr(\Phi)$ are topological group morphisms if the right-hand side is endowed with the compact-open $C^k$-topology for $k\in \N\cup \{\infty\}$ (cf.~\cite[Section 3]{NaS13} for the topological group structure)

Recall however, that since $\Fr(\cO)$ is non-compact, none of these topologies induce a Lie group structure on the group of $C^k$-diffeomorphisms. If we endowed the group with the finer function space topologies which induce a Lie group structure in the smooth case, then the image of $F_k$ is discrete as every point is contained in its own connected component. Due to the compact support conditions enforced by the fine topologies, cf.~\cite{HaS17} as $\Fr(f)\neq\Fr(g)$ differ on a whole non-compact $\GL(n)$-orbit, they differ outside every compact set. 
\end{remark}

\bibliographystyle{new}
\bibliography{orbifolds} 
\end{document}